\documentclass[final,onefignum,onetabnum]{siamart250211}

\usepackage{colortbl}
\definecolor{LightGray}{rgb}{0.83, 0.83, 0.83}

\usepackage{hyperref}       
\hypersetup{colorlinks=true, linkcolor=blue, breaklinks=true, urlcolor=blue}

\usepackage{amssymb,amsmath,amsfonts,mathtools}
\usepackage{stmaryrd}
\usepackage[numbers,sort,compress]{natbib}

\usepackage{enumerate,enumitem}
\renewcommand{\theenumi}{(\roman{enumi})}

\newcommand\oprocendsymbol{\hbox{$\square$}}
\newcommand\oprocend{\relax\ifmmode\else\unskip\hfill\fi\oprocendsymbol}

\newsiamremark{remark}{Remark}
\crefname{remark}{Remark}{Remarks}
\newsiamremark{problem}{Problem}
\crefname{problem}{Problem}{Problems}

\newcommand{\real}{\mathbb{R}}
\newcommand{\complex}{\mathbb{C}}
\newcommand{\bbK}{\mathbb{K}}
\newcommand{\jac}[1]{D\mkern-2.5mu{#1}}
\newcommand{\until}[1]{\{1,\dots, #1\}}
\newcommand{\diag}{\operatorname{diag}}

\newcommand{\norm}[2]{\|#1\|_{#2}}
\newcommand{\WP}[2]{\left\llbracket{#1}, {#2}\right\rrbracket}
\newcommand{\SIP}[2]{\left[{#1}, {#2}\right]}

\newcommand{\setdef}[2]{\{#1 \; | \; #2\}}
\newcommand{\map}[3]{#1: #2 \rightarrow #3}

\newcommand{\realpart}{\operatorname{{Re}}}
\newcommand{\imagpart}{\operatorname{{Im}}}

\newcommand{\barAA}{H}
\newcommand{\Iinfty}{I_{\infty}}

\newcommand{\sign}{\operatorname{sign}}

\newcommand{\subject}{\text{subject to}}

\DeclareSymbolFont{bbold}{U}{bbold}{m}{n}
\DeclareSymbolFontAlphabet{\mathbbold}{bbold}
\newcommand{\vect}[1]{\mathbbold{#1}}

\def\Mili{{Mili\v ci\'c}}

\newcommand{\Id}{\mathbf{I}}
\newcommand{\lognorm}[2]{\mu_{#2}(#1)}
\newcommand{\ymax}{y_{\max}}

\usepackage{stackengine}
\stackMath
\newlength\matfield
\newlength\tmplength
\def\matscale{1.}
\newcommand\dimbox[3]{%
  \setlength\matfield{\matscale\baselineskip}%
  \setbox0=\hbox{\vphantom{X}\smash{#3}}%
  \setlength{\tmplength}{#1\matfield-\ht0-\dp0}%
  \fboxrule=1pt\fboxsep=-\fboxrule\relax%
  \fbox{\makebox[#2\matfield]{\addstackgap[.5\tmplength]{\box0}}}%
}
\newcommand\raiserows[2]{%
   \setlength\matfield{\matscale\baselineskip}%
   \raisebox{#1\matfield}{#2}%
}
\newcommand\matbox[5]{
  \stackunder{\dimbox{#1}{#2}{\ensuremath{#5}}}{\scriptstyle(#3\times #4)}%
}

\title{Regular Pairings for Non-quadratic Lyapunov Functions and Contraction Analysis}

\author{Anton V. Proskurnikov\thanks{Department of Electronics and
    Telecommunications, Politecnico di Torino, Turin, Italy
    (\email{anton.p.1982@ieee.org}). Anton Proskurnikov is supported by the project 2022K8EZBW ``Higher-order interactions in social dynamics with application to monetary networks'', funded by European Union – Next Generation EU within the PRIN 2022 program (D.D. 104 - 02/02/2022 Ministero dell'Universit\'a e della Ricerca). This manuscript reflects only the
authors' views and opinions, and the Ministry for Universities and Research cannot be considered responsible for them.}
    \and Francesco Bullo\thanks{Center
    for Control, Dynamical Systems, and Computation, UC Santa Barbara,
    Santa Barbara, CA 93101 USA (\email{bullo@ucsb.edu}). The work by Francesco Bullo was
    supported in part by AFOSR grant FA9550-22-1-0059. }}

\headers{Regular Pairings for Contraction Analysis}{Proskurnikov and Bullo}

\begin{document}
\maketitle

\begin{abstract}
  Recent studies on stability and contractivity have highlighted the
  importance of semi-inner products, which we refer to as ``pairings'',
  associated with general norms. A pairing is a binary operation that
  relates the derivative of a curve's norm to the radius-vector of the
  curve and its tangent. This relationship, known as the curve norm
  derivative formula, is crucial when using the norm as a Lyapunov
  function. Another important property of the pairing, used in stability
  and contraction criteria, is the so-called Lumer inequality, which
  relates the pairing to the induced logarithmic norm.  We prove that the
  curve norm derivative formula and Lumer's inequality are, in fact,
  equivalent to each other and to several simpler properties. We then
  introduce and characterize \emph{regular pairings} that satisfy all of
  these properties. Our results unify several independent theories of
  pairings (semi-inner products) developed in previous work on functional
  analysis and control theory.  Additionally, we introduce the polyhedral
  max pairing and develop computational tools for polyhedral norms,
  advancing contraction theory in non-Euclidean spaces.
\end{abstract}




\section{Introduction}

\subsection*{Problem description and motivation}

This article is motivated by the broad and growing interest in the
stability and contractivity properties of dynamic systems defined on vector
spaces with norms. In such problems, it is natural to consider Lyapunov
functions defined using the ambient norm.  When computing their Lie
derivative along the system flow, an appropriate binary operation on the
space, which we refer to as a ``pairing'', naturally arises. A pairing is a
generalization of a conventional inner product on the vector space that
satisfies a reduced set of axioms. In summary, this paper focuses on pairings on normed spaces and two key properties they
exhibit: the \emph{curve norm derivative formula} and \emph{Lumer's
property}~\cite{AD-SJ-FB:20o}.

Building on rich and insightful historical developments (reviewed in the next section), recent applications have identified the minimal set of axioms required to perform the Lie product computation and have revealed close relationships between pairings and the induced logarithmic norm on the normed space.
These recent
works have facilitated the systematic treatment of contraction
theory\footnote{A discrete-time dynamical system is \emph{contracting} if
its update map is a contraction in some metric. Analogously, a
continuous-time system is contracting if its flow map is a
contraction. Contraction theory for dynamical
systems~\cite{WL-JJES:98,ZA-EDS:14b,HT-SJC-JJES:21,FB:24-CTDS} provides a set of
concepts and tools for the study and design of continuous and discrete-time
dynamical systems.}  for dynamical systems over vector spaces with general
non-Euclidean norms. In particular, the \emph{sign pairing} and the
\emph{max pairing} have been defined, which are associated, respectively,
with the widely used $\ell_1$ and $\ell_\infty$ norms on $\mathbb{R}^n$.
This emerging framework has recently been applied to control, learning, and
optimization problems; example results based upon the $\ell_1$ and
$\ell_\infty$ norms include contractivity of recurrent neural
networks~\cite{SJ-AD-AVP-FB:21f,AD-AVP-FB:22q}, a non-Euclidean
S-Lemma~\cite{AD-SJ-AVP-FB:23b}, and a monotone operator
theory~\cite{AVP-AD-FB:22a}.

Despite this progress, fundamental open questions remain regarding which
specific pairings satisfy the curve norm derivative formula and Lumer's
property and the relationship between various relevant definitions of
pairings (or semi-inner products) in the literature. There is also a
natural desire to enrich the theory with additional examples of norms
relevant to applications.
Motivated by this analysis, this paper provides a comprehensive, unifying
characterization of pairings satisfying a broad range of equivalent useful
properties. We call such pairings \emph{regular}.
Deferring a detailed historical discussion to Section II, we briefly note that the concept of regular pairings encompasses, as special cases, the directional derivative of a norm--introduced independently by James~\cite{RCJ:47}, \Mili~\cite{PMM:71,PMM:73} and Tapia~\cite{RAT:73} and referred to henceforth as the JMT pairing--and the semi-inner product in the sense of Lumer~\cite{GL:61,GL-RSP:61} and Giles~\cite{JRG:67}, which we term the LG pairing.
Additionally, we define the notion of polyhedral max pairing and develop some computational tools for polyhedral norms.

\subsection*{Applications of pairings and general norms}

The interest for general norms on vector spaces (e.g., the non-differentiable $\ell_1$, $\ell_\infty$ and polyhedral
norms) and computational tools
adapted to them is motivated by applications to networked systems, such as
compartmental systems~\cite{JAJ-CPS:93}, biological transcriptional
systems~\cite{GR-MDB-EDS:10a}, Hopfield neural
networks~\cite{YF-TGK:96,HQ-JP-ZBX:01,SJ-AD-AVP-FB:21f,AD-AVP-FB:21k,AD-AVP-FB:22q},
reaction networks~\cite{MAAR-DA:16,FB-GG:15}, traffic
networks~\cite{SC-MA:15,GC:17,SC:19}, vehicle
platoons~\cite{JM-GR-RS:19}, and coupled
oscillators~\cite{GR-MDB-EDS:13,ZA-EDS:14}.  References on polyhedral
Lyapunov functions and polyhedral set invariance
include~\cite{APM-ESP:86,AB-CB:88,MV-GB:89,GB:91,FB:91,HK-JA-PS:92,FB:95,AP:95,AP:97,FB:99}.
Other examples of non-Euclidean norms that should be mentioned include Barabanov (or extremal) norms~\cite{RT-MM:12},
polynomial norms, defined by the $d$th root of a degree-$d$ homogeneous
polynomial~\cite{AAA-EdK-GH:19}, and canonical dilation-induced homogeneous
norms~\cite{AP:18}.

Recent work on the stability and contractivity of dynamical systems has highlighted several advantages of non-Euclidean norms over standard weighted Euclidean norms. First, polyhedral Lyapunov functions are required~\cite{FB-GG:15} to establish structural stability in biochemical reaction networks, as quadratic Lyapunov functions prove inadequate. Furthermore,~\cite{ZA-EDS:13} presented a class of biochemical models for which stability can be established using the weighted $\ell_1$ norm as a Lyapunov function, whereas weighted $\ell_p$ norms with $p>1$ appear ineffective. Second, the $\ell_1$ and $\ell_\infty$ norms enable scalable and distributed stability and contractivity tests of large-scale networks~\cite{YF-TGK:96,AD-AVP-FB:21k}, a property that linear matrix inequalities for quadratic norms do not offer. Third, the stability properties of monotone dynamical flow networks are readily and generally established~\cite{GC:17} using the $\ell_1$ norm.

\subsection*{Main contributions and paper organization}

Section~\ref{sec:prelim} establishes our notation and reviews preliminary material on pairings and logarithmic norms.  In Section~\ref{sec:reg-pairings}, we present our main results.  We begin in Subsection~\ref{subsec:defs} by introducing the key properties one typically seeks in a pairing—Lumer’s property, the curve-norm derivative formula, and several other equivalent characterizations.  These concepts set the stage for our central \emph{characterization theorem}, proved in Subsection~\ref{subsec:mainthm}, where we demonstrate that all of these conditions are indeed equivalent.  This equivalence motivates our definition of a \emph{regular pairing}: any weak pairing satisfying one (and hence all) of the equivalent properties.  The remainder of Section~\ref{sec:reg-pairings} is devoted to illustrating this concept with important examples: upper JMT pairings, LG pairings, pairings induced by weighted  $\ell_p$ norms, and polyhedral pairings arising from polyhedral norms.

Section~\ref{subsec:motivating-section} reviews the fundamental theorem of contraction theory through the lens of our results on regular pairings, with particular emphasis on the curve-norm derivative formula and Lumer’s inequality.
At the end of this section, we provide an example illustrating that, even in a two-dimensional system, Euclidean norms may be insufficient to establish contraction, whereas a suitable non-Euclidean norm succeeds. Section~\ref{sec:concl} concludes the paper, and all technical proofs are collected in the Appendices.

\section{Preliminaries}\label{sec:prelim}

In this section, we introduce the basic definitions and notation that will be used throughout the text.
We use $\jmath\in\mathbb{C}$ to denote the imaginary unit. The conjugate of a complex number $z\in\mathbb{C}$ is denoted by $\bar z$.

We use $\|\cdot\|_p$, $p\in[1,\infty]$, to denote the conventional $\ell_p$-norm on $\mathbb{R}^n$. The weighted $\ell_p$-norm
is $\|x\|_{p,R}:=\|Rx\|_p$, where $R$ is an invertible square matrix. The corresponding logarithmic norms (see Definition~\ref{def.lognorm}) are denoted
by $\mu_p$ and $\mu_{p,R}$.

\subsection*{Operators and Induced Norms}

Henceforth, $X$ denotes a vector space over the field of scalars $\bbK$
(which is either $\real$ or $\complex$).  We use $\Id_X$ to denote the
identity map on $X$ and denote the zero vector by $\vect{0}$.
If $X$ is equipped by a norm $\|\cdot\|$, we use $\mathcal{B}(X)$ to denote the set of bounded linear operators on $X$.
We will broadly use the induced operator and logarithmic norms of a bounded
operator~\cite{GL:61,GS:24}.
\begin{definition}[\bf The Induced Operator Norm and Logarithmic Norm]
	\label{def.lognorm}
  For an operator $A\in\mathcal{B}(X)$ on a normed space $(X,\|\cdot\|)$, we define
  \begin{itemize}
  \item the induced operator norm
  $
  \|A\|:=\sup_{x\in X:\|x\|=1}\|Ax\|,
  $ and
  \item the induced logarithmic norm, or log norm
  \begin{equation}\label{eq.log-norm}
    \mu(A):=\lim_{h\to 0+}\frac{\|\Id_X+hA\|-1}{h}.
  \end{equation}
  \end{itemize}
\end{definition}

The limit in~\eqref{eq.log-norm} always exists due to the convexity of the norm.
Log norms on $\mathcal{B}(X)$ have practically the same properties as log
norms on square complex matrices (known also as ``matrix measures'') that
are surveyed in~\cite{CAD-MV:1975,GS:24,FB:24-CTDS}.

\subsection*{Pairings (Semi-Inner Products)}

In this work, the term \emph{pairing} refers to a binary operation\footnote{A \emph{binary operation} on $X$ is a function $\SIP{\cdot}{\cdot}:X\times
X\to\bbK$.}  on a vector space that intentionally satisfies only a subset of the inner product axioms, being
homogeneous in its first argument, positive on the diagonal, and satisfies a generalized form of the Cauchy–Schwarz inequality, while not necessarily being linear or symmetric. While the term ``semi-inner product'' (SIP) is frequently used in the literature to denote such structures, we avoid it here due to its ambiguity -- a point that will become evident in the historical discussion that follows.

In this work, we focus on three principal classes of binary operations on normed vector spaces: LG pairings, JMT pairings, and weak pairings (WP).
\begin{definition}[\bf LG Pairing, or the Lumer-Giles SIP]
  A semi-inner product in the sense of Lumer and Giles, or briefly
  \emph{LG} pairing, on the vector space $X$ is a binary operation
  $\SIP{\cdot}{\cdot}$ satisfying the three axioms:
  \begin{enumerate}
  \item Linearity in the first argument:
    \begin{gather}
    \SIP{x_1+x_2}{y}=\SIP{x_1}{y}+\SIP{x_2}{y}\qquad\forall x_1,x_2,y\in X;\label{eq.additive}\\
    \SIP{\lambda x}{y}=\lambda\SIP{x}{y}\qquad\forall \lambda\in\bbK\;\;\forall x,y\in X;\label{eq.hom1}
    \end{gather}
  \item Positive definiteness:
    \begin{equation}\label{eq.posdef-sip}
      \SIP{x}{x}>0\qquad\forall x\in X\setminus\{\vect{0}\};
    \end{equation}
  \item The Cauchy-Schwarz inequality:
    \begin{equation}\label{eq.cauchy-sip}
      \SIP{x}{y}\leq \sqrt{\SIP{x}{x}\SIP{y}{y}}\qquad\forall x,y\in X;
    \end{equation}
  \item Complex-conjugate homogeneity in the second argument:
    \begin{equation}\label{eq.hom2}
      \SIP{x}{\lambda y}=\bar{\lambda}\SIP{x}{y}\qquad\forall\lambda\in\bbK\;\;\forall x,y\in X.
    \end{equation}
  \end{enumerate}
\end{definition}

In his first works Lumer~\cite{GL:61,GL-RSP:61} defined a SIP as a binary
operation satisfying the axioms \emph{(i)-(iii)}; his definition resembles
the usual inner however, the (conjugate) symmetry axiom
$\SIP{x}{y}=\overline{\SIP{y}{x}}$ is discarded.  It can be easily shown
that, similar to usual inner products, such a binary operation naturally
induces a norm
\begin{equation}\label{eq.norm-induced}
  \|x\|:=\sqrt{\SIP{x}{x}},
\end{equation}
and, on the other hand, every norm on $X$ can be obtained in this way.
Axiom \emph{(iv)} (which is often convenient in applications) was later
introduced by Giles~\cite{JRG:67} who showed that every norm could in fact
be represented by a SIP satisfying \emph{(i)-(iv)}.

\begin{definition}[{\bf JMT Pairings, or the James-\Mili-Tapia SIPs}]
  Consider a vector space $X$ equipped with a norm $\|\cdot\|$. The
  \emph{upper} and \emph{lower semi-inner products} in the sense of James,\footnote{James~\cite{RCJ:47} introduced the left and right Gateaux derivatives $N_+(y,x)$ and $N_-(y,x)$ of a general norm. Adopting his notation, $\SIP{x}{y}_{+}=N_{+}(y,x)\|y\|$ and $\SIP{x}{y}_{-}=N_{-}(y,x)\|y\|$.}
  Mili\v ci\'c\footnote{Mili\v ci\'c~\cite{PMM:71} introduced another
  generalization of the inner product, which is sometimes also called the
  semi-inner product in the Mili\v ci\'c sense (M-SIP)~\cite[Chapter~4,
    Definition~8]{SSD:04} and is defined as the average of the upper and
  the lower JMT pairings~\eqref{eq.JMT-SIP}. This binary operation, in general,
  does not belong to the class of weak pairings defined below and is beyond the scope of our work.} and Tapia (termed JMT pairings)
  induced by $\|\cdot\|$ are the \emph{real-valued} binary operations
  $\SIP{x}{y}_+$, $\SIP{x}{y}_-$ defined by
  \begin{equation}\label{eq.JMT-SIP}
    \SIP{x}{y}_+:=\|y\|\lim_{h\to 0+}\frac{\|y+hx\|-\|y\|}{h},\quad \SIP{x}{y}_-:=\|y\|\lim_{h\to 0-}\frac{\|y+hx\|-\|y\|}{h}.
  \end{equation}
\end{definition}

Using the convexity of the norm, it can be shown~\cite{RCJ:47,PMM:71} that both
limits in~\eqref{eq.JMT-SIP} exist, furthermore,
$\SIP{x}{y}_-\leq\SIP{x}{y}_+=-\SIP{-x}{y}_-$ for all $x,y\in X$.

\begin{remark}[\bf Differences Between the LG and JMT Pairings]
  It is useful to review the principal differences between the LG and JMT
  pairings~\cite{SSD:04}. In complex spaces, LG pairings are
  complex-valued, whereas JMT pairings can only be real-valued. Whereas LG
  pairings can be defined on an arbitrary vector space and induces a norm
  on it in accordance with~\eqref{eq.norm-induced}, the JMT pairings can be
  defined only on a normed space, providing, however, the compatibility
  between the SIP and the norm
  \[
  \|x\|^2=\SIP{x}{x}_+=\SIP{x}{x}_-.
  \]
  The upper and lower JMT pairings obey~\eqref{eq.hom1}-\eqref{eq.hom2} yet do not satisfy~\eqref{eq.additive}, being, instead, subadditive and superadditive in the first argument respectively
  \begin{align}
    \SIP{x_1+x_2}{y}_+ &\leq \SIP{x_1}{y}_+ + \SIP{x_2}{y}_+\label{eq.subadd-jmt}\\
    \SIP{x_1+x_2}{y}_- &\geq \SIP{x_1}{y}_- + \SIP{x_2}{y}_-.
  \end{align}
  Finally, the JMT pairings are not homogeneous: in fact,
  equations~\eqref{eq.hom1} and~\eqref{eq.hom2} hold only for $\lambda\geq
  0$. However, it can be shown that for any $\lambda\in\bbK$ one has
  \begin{equation*}
    \SIP{\lambda x}{\lambda y}_{+}=|\lambda|^2\SIP{x}{y}_{+},\; \SIP{\lambda x}{\lambda
      y}_{-}=|\lambda|^2\SIP{x}{y}_{-}\qquad\forall x,y\in X.
  \end{equation*}
  This concludes our comparison. \oprocend
\end{remark}

Another definition of a pairing, which enjoys particular properties of both LG and
JMT pairings, was proposed in~\cite{AD-SJ-FB:20o}, which paper has also
coined a term \emph{weak pairing} (WP) for this class of operations.
Formally, the definition from~\cite{AD-SJ-FB:20o} was confined to $\real^n$
but can be trivially extended to an arbitrary space.
\begin{definition}[{\bf Weak Pairing}]
  A real-valued binary operation $\WP{\cdot}{\cdot}$ on $X$ is said to be a
  weak pairing (WP) if it satisfies the following axioms\footnote{Notice
  that in~\cite{AD-SJ-FB:20o} (considering $X=\real^n$) an additional
  assumption on the continuity in the first argument was imposed. Lemma~\ref{lem.wp-induces-norm} shows that this assumption was in fact
  superfluous as all pairings on a finite-dimensional space are Lipschitz continuous in the first argument.}
\begin{enumerate}
\item\label{WP1}(Subadditivity in the first argument)
\begin{equation}\label{eq.subadd}
\WP{x_1+x_2}{y} \leq \WP{x_1}{y} + \WP{x_2}{y}\qquad\forall x_1,x_2,y \in X;
\end{equation}
\item\label{WP3}(Weak homogeneity)
\begin{gather}
\WP{\alpha x}{y} = \WP{x}{\alpha y} = \alpha\WP{x}{y}\qquad\forall x,y\in X,\qquad\forall \alpha\geq 0\label{eq.weak-hom1}\\
\WP{-x}{-y} = \WP{x}{y}\qquad\forall x,y\in X.\label{eq.weak-hom2}
\end{gather}
\item\label{WP4}(Positive definiteness) $\WP{x}{x} > 0$, for all $x \neq \vect{0}$
\item\label{WP5}(Cauchy-Schwarz inequality)
\begin{equation}\label{eq.cauchy-wp}
\WP{x}{y}\leq \sqrt{\WP{x}{x}\WP{y}{y}}\qquad\forall x,y\in X;
\end{equation}
\end{enumerate}
\end{definition}

The upper JMT pairing and the LG pairing are special cases of a weak pairing, which inherits properties from both. Like the JMT pairing, the weak pairing is real-valued, subadditive, and satisfies a weak form of homogeneity. Simultaneously, similar to the LG pairing, the weak pairing is defined on a general vector space and automatically induces a norm, as demonstrated by the following.
\begin{lemma}[{\bf Weak Pairings and Norms}]\label{lem.wp-induces-norm}
  For every WP $\WP{\cdot}{\cdot}$ on $X$, the function
  $\|x\|:=\sqrt{\WP{x}{x}}$ is a norm. In this norm, the
  function $\WP{\cdot}{y}:X\to\real$ is Lipschitz continuous for every $y\in X$, and its Lipschitz constant is $\|y\|$, that is,
  \begin{equation}\label{eq.wp-lip1}
    |\WP{x_1}{y}-\WP{x_2}{y}|\leq \|x_1-x_2\|\,\|y\|\quad\forall x_1,x_2,y\in X.
  \end{equation}
\end{lemma}
\begin{proof}
  The first statement is proved similarly
  to~\cite[Theorem~16]{AD-SJ-FB:20o}. To prove the second statement, notice
  that for all vectors $x_1,x_2,y\in X$ the subadditivity, the
  Cauchy-Schwarz inequality and the definition of the norm imply that
  \[
  \WP{x_1}{y}\overset{\eqref{eq.subadd}}{\leq} \WP{x_1-x_2}{y}+\WP{x_2}{y}\overset{\eqref{eq.cauchy-wp}}{\leq} \|x_1-x_2\|\,\|y\|+\WP{x_2}{y}.
  \]
  which, in turn, implies that $\WP{x_1}{y}-\WP{x_2}{y}\leq \|x_1-x_2\|\,\|y\|$. Swapping $x_1$ and $x_2$, one proves similarly that
  $\WP{x_2}{y}-\WP{x_1}{y}\leq \|x_1-x_2\|\,\|y\|$, which entails~\eqref{eq.wp-lip1}.
\end{proof}

\subsection*{Historical discussion on pairings and applications}

Motivated by efforts to extend the well-established theory of Hilbert spaces to broader classes of normed spaces~\cite{AVB:81,PRH:82}, the theory of \emph{semi-inner} product (SIP) spaces was developed. In its broadest sense, the term SIP refers to a binary operation on a vector space that satisfies only a subset of the axioms of an inner product while omitting others.  In some functional analysis
monographs (e.g., in~\cite{JBC:90}) the term SIP is used to denote a binary
operation that adheres to all the axioms of an inner product except for
strict positivity; such an operation is linked with a seminorm rather than
a norm. In this paper, we explore alternative types of SIPs that are positive
definite and linked to norms but are characterized by linearity (or even
sublinearity) in only one argument and the absence of symmetry. To avoid
ambiguity, we refer such SIPs to as \emph{pairings}.
The three key classes of pairings discussed in this paper were introduced independently in the literature, each emerging from distinct lines of research.

The \textbf{first} direction of research, leading to the JMT pairings, was initiated in the seminal work by James~\cite{RCJ:47} on
orthogonality in normed spaces, followed by contributions by
\Mili~\cite{PMM:71,PMM:73} and Tapia~\cite{RAT:73}, on characterization of
inner product spaces.  Defined on a real or a complex normed space, JMT pairings
are real-valued, continuous, semi-linear with respect to one argument (the
direction of the derivative), and positive definite, however, they are not
conjugate-symmetric. Deimling~\cite[Ch.~3]{KD:85} established several
important properties of JMT pairings, including the relationship between
them and the derivative of the norm of a curve in a Banach space;
in~\cite{AD-SJ-FB:20o}, JMT pairings were also called ``Deimling's
pairings''. A further review of the JMT pairing theory is available
in~\cite{SSD:04}. The formalism of JMT pairings enables the introduction of the one-sided
Lipschitz constant of a mapping \cite{GS:24} and establishes a
straightforward relationship with the logarithmic norm of the mapping's
Jacobian, if it exists. This theory extends naturally to more general weak
pairings~\cite{AD-SJ-FB:20o} and non-differentiable locally Lipschitz
maps~\cite{AD-AVP-FB:22q}. Aminzare and Sontag~\cite{ZA-EDS:14,ZA-EDS:14b,ZA-EDS:13} were the first to establish
important connections between the logarithmic norms, JMT pairings and contraction theory for dynamical systems;
this theory has been further developed in~\cite{AD-SJ-FB:20o,FB:24-CTDS}.
JMT pairings also illuminate the properties of accretive and dissipative
operations in general Banach spaces. The accretivity property can be
expressed as a standard monotonicity condition, with the inner product
replaced by a JMT pairing~\cite{KD:85}.  This theory has been further
developed in~\cite{AD-SJ-AVP-FB:23b}, where the concept of monotonicity in
a general non-Euclidean norm is introduced using weak pairings.

On a parallel \textbf{second} line of research,
Lumer~\cite{GL:61,GL-RSP:61} introduced an alternative definition of a SIP
with the aim of extending certain results and concepts from operator
theory--such as dissipative $C_0$-semigroups, Hermitian operators and
associated Hermitian forms--to Banach spaces that lack an inner product. A
SIP, as defined by Lumer, can be defined on both real and complex normed
spaces, and, unlike the JMT pairings, is \emph{linear} in its first
argument. This SIP also lacks symmetry (except in case where it coincides
with a conventional inner product). Additionally, Lumer's definition
permits non-homogeneity in its second argument. Giles~\cite{JRG:67} further
advanced Lumer's concept and proved that each norm is generated by a SIP
that possesses homogeneity in both arguments. By incorporating the convenient
homogeneity axiom into Lumer's definition of a SIP, the notion of the
LG pairing presented above emerges. Applications of LG pairings include elegant characterizations of isometric
operators on normed space~\cite{DK-PR:70,PW-TS:21} and a refined theory of
semi-polarity in non-Euclidean norms~\cite{AGH-ZL-MS:15}. Inspired by the
applications of Hilbert spaces in machine learning, signal processing and
numerical analysis, the concepts of reproducing kernel Banach
spaces~\cite{HZ-JZ:10} and Riesz bases in Banach
spaces~\cite{HZ-JZ:11} have been introduced.

The \textbf{third} line of research, developed recently by Davydov et
al.~\cite{AD-SJ-FB:20o}, introduced the concept of a \emph{weak pairing},
an axiomatic definition with weaker assumptions than the LG pairing
notion. The \emph{sign pairing} (naturally associated with the $\ell_1$
norm) and the \emph{max pairing} (generating the $\ell_{\infty}$ norm) are
examples of weak pairings. In fact, every LG and JMT pairing on a real
normed space is a special case of a weak pairing. Aimed at developing a
unified contraction theory in normed spaces,~\cite{AD-SJ-FB:20o} reveals
the minimum critical properties required from a weak pairing to establish
Lyapunov-based contractivity criteria, namely, Lumer's property, the curve
norm derivative formula, and the dominance by the upper JMT pairing, or the
``Deimling inequality''~\cite{AD-SJ-FB:20o}. In the next section, we establish that these properties are, in fact, equivalent; a weak pairing satisfying any of them is termed \emph{regular}. This equivalence is formalized in Theorem~\ref{thm.main}, which provides several mutually equivalent characterizations of regularity. Notably, both the Lumer–Giles (LG) and upper JMT pairings satisfy these conditions, and thus qualify as regular pairings.  In this sense, the paper brings together and systematizes three previously distinct lines of research.

\section{Regular Pairings and Their Characterization}
\label{sec:reg-pairings}

In this section, we will show that some useful properties established for
LG and JMT pairings~\cite{SSD:04} and some special WPs~\cite{AD-SJ-FB:20o} are
mutually equivalent and are enjoyed by a broad class of WP that we call
\emph{regular pairings}.  The key Theorem~\ref{thm.main} below generalizes,
as will be shown, a number of results on semi-inner products available in
the literature.

Henceforth, $(X,\|\cdot\|)$ is a normed space whose norm is associated with a WP $\WP{\cdot}{\cdot}$.

\subsection{Technical Definitions: Key Properties of a Pairing}\label{subsec:defs}

We start with auxiliary definitions and introduce key properties desired
from a pairing~\cite{AD-SJ-FB:20o}.
\begin{definition}[{\bf Lumer's Property}]
  Following~\cite{AD-SJ-FB:20o}, we say that \emph{Lumer's property} holds
  for the WP if, for each linear bounded operator $A:X\to X$,
  \begin{equation}\label{eq.lumer}
    \mu(A)=\sup_{x\,:\,\|x\|=1}\WP{Ax}{x}.
  \end{equation}
\end{definition}

Recalling that a WP is weakly homogeneous~\eqref{eq.weak-hom1}, by
  substituting $x=z/\|z\|$, where $z\ne 0$, Lumer's property can be
  reformulated as follows:
  \begin{equation}\label{eq.lumer-a}
    \mu(A)=\sup_{z\ne 0}\frac{\WP{Az}{z}}{\|z\|^2}.
  \end{equation}

\begin{definition}[\bf Curve Norm Derivative]
  The WP enjoys the \emph{curve norm derivative} property if the norm
  $\|x(\cdot)\|$ of every differentiable curve $x:(a,b)\to X$ (where
  $a,b\in\real,\,a<b$) satisfies the following equality
  \begin{equation}\label{eq.curve-norm}
    \frac{1}{2}\frac{d}{dt}(\|x(t)\|^2)=\|x(t)\|\frac{d}{dt}\|x(t)\|=\WP{\dot x(t)}{x(t)}.
  \end{equation}
  at any point $t\in(a,b)$ where the left-hand side is
  defined\footnote{According to Lemma~\ref{lem.curve-norm-deriv}, the norm
  of a differentiable curve is differentiable at almost any point.}.
\end{definition}

The importance of the curve norm derivative formula will be further
clarified in Section~\ref{subsec:motivating-section}. In stability and contraction
analysis for dynamical systems over normed spaces, this formula enables the
use of the squared norm $V(x)=\frac{1}{2}\|x\|^2$ as a Lyapunov function, explicitly
computing its Lie derivative.

\begin{definition}[\bf Partial Linearity]
  The WP is \emph{partially linear} in its first argument if
  $\WP{x+ay}{y}=\WP{x}{y}+a\|y\|^2$ for all $x\in X,y\in X,a\in\real$.
\end{definition}

\begin{definition}[\bf Straight Angle Property]
  The WP satisfies the \emph{straight angle property} if
  $\WP{-x}{x}=-\|x\|^2$ for all $x\in X$.
\end{definition}

The term ``straight angle property'' is inspired by an analogy with
inner-product spaces: for an inner product $(\cdot,\cdot)$, the angle
between vectors $x,y\ne\vect{0}$ is defined as
\[
\angle(x,y)=\arccos\frac{\realpart(x,y)}{\|x\|\,\|y\|}.
\]
Generalizing this definition to WPs, one can reformulate the relation
$\WP{-x}{x}=-\|x\|^2$ as the requirement $\angle(-x,x)=\pi$ whenever $x\ne\vect{0}$. Notice
that, in general, WP may fail to satisfy this condition and, moreover, may
attain only \emph{nonnegative} values. This is illustrated, e.g., by the
following WP on $\real^n$, compatible with the Euclidean norm:
\begin{equation*}
\WP{x}{y}=\sum_{i=1}^n|x_iy_i|.
\end{equation*}

\subsection{Main result: The Characterization Theorem and Regular Pairings}\label{subsec:mainthm}

We now formulate the main result of this section, which will be proved in Appendix~\ref{app.proof-thm}.

\begin{theorem}[\bf The Characterization Theorem]\label{thm.main}
Let $(X,\|\cdot\|)$ be a normed space and $\WP{\cdot}{\cdot}$ be a WP
compatible with the norm: $\|x\|^2=\WP{x}{x}$.  Then, the
following statements are equivalent:
\begin{enumerate}
  \item\label{thm1:straight} the WP obeys the straight angle property;
  \item\label{thm1:linear} the WP is partially linear in its first argument;
  \item\label{thm1:dominated} the WP is related to the JMT pairings~\eqref{eq.JMT-SIP}, corresponding to the norm $\|\cdot\|$, by the following inequalities
  \begin{equation}\label{eq.dominated}
  \SIP{x}{y}_-\leq -\WP{-x}{y}\leq\WP{x}{y}\leq \SIP{x}{y}_+\qquad\forall x,y\in X;
  \end{equation}
  \item\label{thm1:lumer1} for each bounded linear operator $A$, the ``one-sided'' Lumer inequality holds
  \begin{equation}\label{eq.lumer1}
    \WP{Ax}{x}\leq\mu(A)\qquad\forall x\in X: \|x\|=1;
  \end{equation}
    \item\label{thm1:curve1} the curve norm derivative formula~\eqref{eq.curve-norm} holds for all differentiable curves $x:(a,b)\to X$;
  \item\label{thm1:curve2} the equality~\eqref{eq.curve-norm} holds for affine functions $x(t)=tv+x_0$ defined on $(a,b)=\real$ (here $v,x_0\in X$), that is,
      \begin{equation}\label{eq.curve-norm1}
      \frac{1}{2}\frac{d}{dt}\|x_0+vt\|^2=\|x_0+vt\|\left(\frac{d}{dt}\|x_0+vt\|\right)=\WP{v}{x_0+vt}
      \end{equation}
      at each point $t\in(a,b)$ where the left-hand side is well defined;
  \item\label{thm1:lumer} the WP enjoys the Lumer property.
\end{enumerate}
\end{theorem}

\begin{remark}[\bf Absolutely Continuous Curves]
  In $X=\bbK^n$, a differentiable curve
  $x(\cdot)$ in the statement~\emph{\ref{thm1:curve1}} can be replaced by an
  absolutely continuous\footnote{We call $x:(a,b)\to\complex^n$ absolutely
  continuous if all coordinate functions $x_i(\cdot),i=1,\ldots,n$ are
  absolutely continuous on $(a,b)$ in the standard
  sense~\cite[Definition~7.17]{WR:87}.}  (e.g., locally Lipschitz)
  curve. In this situation, $\|x(\cdot)\|$ is also absolutely continuous
  and~\eqref{eq.curve-norm} holds at each point $t$ where both the
  left-hand side and the right-hand side are well defined (this hold at
  almost any $t$~\cite[Theorem~7.18]{WR:87}). The proof remains same. \oprocend
\end{remark}
\begin{remark}[\bf ``Deimling Inequality''\cite{AD-SJ-FB:20o}]
  The equivalence $\ref{thm1:lumer}\Longleftrightarrow\ref{thm1:dominated}$
  was essentially established in~\cite[Theorem~18]{AD-SJ-FB:20o} for the case
  $X=\real^n$; the rightmost inequality in~\eqref{eq.dominated} is referred in~\cite{AD-SJ-FB:20o} as the ``Deimling inequality''.
  The proof given in~\cite{AD-SJ-FB:20o} generalizes to
  a general normed space $X$ and operator $A\in\mathcal{B}(X)$ whose resolvent $(\Id-hA)^{-1}\in\mathcal{B}(X)$ exists for $h$ being sufficiently small.  If the normed space $X$ is Banach, the latter condition holds as $(\Id-hA)^{-1}=\sum_{j=0}^{\infty}h^jA^j$, which series converges
  when $|h|<\|A\|^{-1}$. Theorem~\ref{thm.main}, however, \emph{does not} rely on the completeness of $X$.\oprocend
\end{remark}

We codify the numerous equivalences in Theorem~\ref{thm.main} via the following definition.

\begin{definition}[\bf Regular Pairing]\label{def.reg-pair}
A \textbf{regular pairing} on the space $X$ is a weak pairing that enjoys the properties~\ref{thm1:straight}-\ref{thm1:lumer} from
  Theorem~\ref{thm.main}.
\end{definition}

Regular pairings -- though not referred to by this name -- have been implicitly used in several prior works to derive contraction criteria for general nonlinear systems~\cite{AD-SJ-FB:20o} and nonlinear neural networks~\cite{AD-AVP-FB:21k,AD-AVP-FB:22q}, robustness properties of implicit machine learning models~\cite{SJ-AD-AVP-FB:21f}, and monotonicity of nonlinear operators in general normed spaces~\cite{AD-SJ-AVP-FB:23b}. However, being primarily focused on weighted $\ell_1$ and $\ell_{\infty}$ norms, these works, however, relied on three key properties of a weak pairing: Lumer's property, the curve norm derivative formula, and the rightmost inequality in~\eqref{eq.dominated} (referred to in~\cite{AD-SJ-FB:20o} as the ``Deimling inequality'').
However, it has not been recognized that these properties are, in fact, equivalent. Their \emph{ad hoc} verification is not entirely straightforward -- even for the max-pairing and sign-pairing associated with the $\ell_1$ and $\ell_{\infty}$ norms~\cite{AD-SJ-FB:20o}. Theorem 3.5 makes it possible to bypass the direct verification of the three key properties by reducing them to much simpler conditions, such as the straight angle property. The theorem entails the regularity of both the LG pairing and the upper JMT pairing, offering a unified explanation of their properties.

The remainder of this section outlines several significant classes of regular pairings, underscoring that Theorem~\ref{thm.main} unifies and extends various results in the literature. An application in control theory is presented in Section~\ref{subsec:motivating-section}, where the contraction criterion (Theorem~\ref{thm:general}) is formulated—restating the earlier result~\cite[Theorem 1]{AD-SJ-FB:20o}; it is demonstrated that the use of non-Euclidean norms, based on regular pairings, can offer significant advantages even for low-dimensional nonlinear systems.

\subsection{Upper JMT Pairings as Regular Pairings}

We begin by observing that the upper JMT pairing induced by a norm constitutes an instance of a regular pairing.
\begin{corollary}[\bf Upper JMT Pairings are Regular Pairings]\label{cor.JMT-SIP-lumer}
  For every normed space $(X,\|\cdot\|)$, the upper JMT pairing
  $\SIP{\cdot}{\cdot}_+$ associated with the norm is a regular pairing on
  $X$. In particular, the upper JMT pairing enjoys the Lumer property.
\end{corollary}
\begin{proof}
  It suffices to notice that $\SIP{\cdot}{\cdot}_+$ obeys~\eqref{eq.dominated} (where the leftmost and the rightmost inequalities turn into equalities).
\end{proof}

Remarkably, there exists a broad class of normed spaces that admit no other regular pairings -- namely, those whose norms are Gateaux differentiable~\cite{RCJ:47,JRG:67,RAT:73}.
\begin{definition}[\bf Gateaux Differentiable Normed Space]\label{def.gateaux}
  A normed space $(X,\|\cdot\|)$ is \emph{Gateaux
  differentiable} if the derivative
  $\biggl.\frac{d}{dt}\biggr|_{t=0}\|y+tx\|$ exists for each non-zero point $y\ne\vect{0}$ and each direction $x\in X$.
  Notice that in this situation $\WP{x}{y}_-=\WP{x}{y}_+$, for all $x,y\in X$.
\end{definition}

An inner product space is always Gateaux differentiable; other examples
include $\ell_p$ and functions $L_p$ spaces for $1<p<\infty$~\cite{JRG:67}.

\begin{lemma}[\bf Unique Regular Pairing for Gateaux Differentiable Norms]\label{thm.unique}
  A norm $\|\cdot\|$ on $X$ is induced by \emph{only one} regular pairing
  if and only if the space $(X,\|\cdot\|)$ is Gateaux differentiable. In
  this case, the only regular pairing compatible with the norm is
  $\WP{\cdot}{\cdot}:=\SIP{\cdot}{\cdot}_+=\SIP{\cdot}{\cdot}_-$ and the
  only LG pairing compatible with the norm is $\WP{\cdot}{\cdot}$ in the real
  case and~\eqref{eq.wp-to-sip} in the complex case.

  In a finite dimensional normed space $(\bbK^n,\|\cdot\|)$, the regular WP is ``almost'' unique in the sense that
  $\WP{\cdot}{y}\equiv\SIP{\cdot}{y}_+\equiv\SIP{\cdot}{y}_-$ for almost all $y\in X$.
\end{lemma}
The proof of Lemma~\ref{thm.unique} is given in Appendix~\ref{app.proof-thm-unique}.

\begin{remark}\label{rem.discont}
  Retracing the arguments from~\cite[Theorem~3]{JRG:67}, it can be shown
  that all regular pairings on non-differentiable (in the Gateaux sense) normed spaces are discontinuous
  in the second argument. Moreover, for some $x,y\in X$ the function $t\mapsto\WP{x}{x+ty}$ is discontinuous at $t=0$. \oprocend
\end{remark}

\begin{remark}
  It should be noticed that Lemma~\ref{thm.unique} is not valid for
  irregular weak pairings. For instance, the Euclidean norm on $\real^n$ is
  induced by pairings
  \[
  \WP{x}{y}_{\alpha}=\alpha\sum\nolimits_{i=1}^nx_iy_i+(1-\alpha)\sum\nolimits_{i=1}^n|x_iy_i|,\quad\alpha\in [0,1],
  \]
  among which the only regular WP corresponds to $\alpha=1$, whereas the others fail to satisfy the straight angle property. \oprocend
\end{remark}

\subsection{Lumer-Giles (LG) Pairings}\label{subsec:lg-jmt-regularity}

In this section, we formulate important corollaries of Theorem~\ref{thm.main} that are concerned with the LG pairings, which are special cases of a regular pairing and, as discussed above, serve as natural extensions of inner products in many applications~\cite{HZ-JZ:10,HZ-JZ:11}.

\begin{corollary}[\bf Real parts of LG Pairings are Regular Pairings]\label{cor.lg-sip}
  The real part\footnote{If the space $X$ is real, the real part can,
  obviously, omitted.} $\WP{x}{y}=\realpart\SIP{x}{y}$ of an LG pairing
  on the space $X$ is a regular pairing on $X$. In particular, this pairing
\begin{enumerate}
\item\label{cor1:dominated} obeys the inequality $\SIP{x}{y}_-\leq \realpart\SIP{x}{y}\leq\SIP{x}{y}_+$ for all $x,y\in X$, where the JMTPs $\SIP{\cdot}{\cdot}_{\pm}$ correspond to the induced norm~\eqref{eq.norm-induced}.
\item\label{cor1:lumer} enjoys the Lumer property~\eqref{eq.lumer}, and
\item\label{cor1:curve} possesses also the curve norm property~\eqref{eq.curve-norm}.
\end{enumerate}
\end{corollary}
\begin{proof}
The proof is straightforward, because $\SIP{-x}{x}=-\SIP{x}{x}=-\|x\|^2$, and hence the straight angle condition holds.
\end{proof}

Notice that statement~\emph{\ref{cor1:dominated}} of Corollary~\ref{cor.lg-sip} is known in the literature, see, e.g.,~\cite[Lemma~2.2]{SSD-JJK:00} and~\cite[Chapter~X, Lemma~5]{SSD:04}. However, we give a direct proof, which does not exploit duality mappings~\cite{SSD-JJK:00}. Statement~\emph{\ref{cor1:lumer}} is the well-known result by Lumer result~\cite[Lemma~12]{GL:61}, relating the log norm and the numerical range of an operator.
Statement~\emph{\ref{cor1:curve}}, however, does not seem easily available in the literature.

A natural question arises: which regular pairings can be represented by
real parts of LG pairings? The answer is given by the following lemma proven in Appendix~\ref{app.proof-lem-wp-sip}.
\begin{lemma}[\bf Regular Pairings Arising from LG Pairings]\label{lem.wp-and-sip}
  For a regular pairing $\WP{\cdot}{\cdot}$ on a real vector space $X$, two
  conditions are equivalent:
  \begin{enumerate}
  \item\label{lem1:real}$\WP{\cdot}{\cdot}$ is a LG pairing on $X$;
  \item\label{lem1:hom} $\WP{-x}{y}=-\WP{x}{y}$ for all $x,y\in X$.
  \end{enumerate}
  For a WP on a complex vector space $X$, the following conditions are
  equivalent:
  \begin{enumerate}\setcounter{enumi}{3}
  \item\label{lem1:real+}$\WP{\cdot}{\cdot}=\realpart\SIP{\cdot}{\cdot}$ is a real part of some LG
    pairing on $X$;
  \item\label{lem1:hom+} $\WP{-x}{y}=-\WP{x}{y}$, $\WP{\jmath x}{x}=0$, and
    $\WP{x}{\lambda y}=\WP{\bar\lambda x}{y}$ for all $x,y\in X$ and
    $\lambda\in\mathbb{C}$ such that $\|\lambda\|=1$.
  \end{enumerate}
  The LG pairing $\SIP{\cdot}{\cdot}$ satisfying~$\ref{lem1:real+}$ is unique.
\end{lemma}



\subsection{Regular Pairings Associated to Weighted $\ell_p$ Norms on $\mathbb{R}^n$}\label{subsec:l_p}

In this section, we focus on regular pairings in finite-dimensional spaces,
extending the theory from \cite{AD-SJ-FB:20o} to highlight the importance
of regular pairings for analysis of dynamical systems. Formally, all norms on a finite-dimensional vector space are equivalent;
however, the constant factors establishing this equivalence increase with
the space's dimension.\footnote{For example, it is known that
$\norm{x}{1}\leq\sqrt{n}\norm{x}{2}\leq{n}\norm{x}{\infty}$ for each
$x\in\real^n$.}  Consequently, estimating a solution with one norm can
become overly conservative when converted to another norm. Therefore, while
positive definite quadratic forms (equivalently, squared weighted $\ell_2$
norms, see below) are the most typical Lyapunov functions, many
applications require tight estimates of system solutions in non-Euclidean
norms. For instance, a long-standing challenge in machine learning is
obtaining certifiable robustness bounds for deep neural networks against
adversarial
perturbations~\cite{CZ-WZ-IS-JB-DE-IG-RF:13,MR-RW-IRM:21,MF-MM-GJP:20},
which are typically measured in the $\ell_{\infty}$
norm~\cite{IJG-JS-CZ:15}. Therefore, the input/output Lipschitz constants
of neural networks in $\ell_{\infty}$ norm need to be tightly estimated.
Other applications, leading to stability and contraction in weighted
$\ell_1$ and $\ell_{\infty}$ norms are surveyed
in~\cite{SJ-AD-AVP-FB:21f,AD-AVP-FB:22q,AD-SJ-AVP-FB:23b}. Among them are
tight estimates for iterative algorithms seeking fixed points or monotone
operators, and efficiently computable contraction rates of recurrent neural
networks, to name a few.

Table~\ref{table:equivalences}, borrowed from \cite{AD-SJ-FB:20o},
summarizes the standard weak pairings compatible with weighted $\ell_p$
norms and the corresponding induced logarithmic norms.
\begin{table*}[h]\centering
  \normalsize
    \resizebox{1\textwidth}{!}
            {\begin{tabular}{%
	p{0.12\linewidth}%
	p{0.39\linewidth}%
	p{0.47\linewidth}%
      }
      Norm & Weak Pairing
      & Logarithmic norm
      \\
      \hline
      \rowcolor{LightGray}    &&  \\[-2ex]
      \rowcolor{LightGray}
      $ \begin{aligned}
	\norm{x}{2,P^{1/2}}
      \end{aligned}$
      &
      $\begin{aligned}
	\WP{x}{y}_{2, P^{1/2}} &= x^\top Py
      \end{aligned}$
      &
      $\begin{aligned}
	&\mu_{2,P^{1/2}}(A) = \tfrac{1}{2}\lambda_{\max}(PAP^{-1} + A^\top)\\
& = \tfrac{1}{2}\lambda_{\max}(P^{1/2}AP^{-1/2} + P^{-1/2}A^\top P^{1/2})\\
	 &= \max_{\|x\|_{2,P^{1/2}}=1} x^\top PAx
      \end{aligned}
      $
      \\[2ex]
      &&  \\[-2ex]
      $ \begin{aligned}
	&\norm{x}{p}
      \end{aligned}$
      &
      $\begin{aligned}
	\WP{x}{y}_{p} &= \|y\|_{p}^{2-p}(y \circ |y|^{p-2})^\top x
      \end{aligned}$
      &
      $\begin{aligned}
	\mu_{p}(A) &= \max_{\|x\|_{p} = 1} (x \circ |x|^{p-2})^\top Ax
      \end{aligned}
      $
      \\
      \rowcolor{LightGray}    &&  \\[-2ex]
      \rowcolor{LightGray}
      $ \begin{aligned}
	\norm{x}{1}
      \end{aligned}$
      &
      $\begin{aligned}
	\WP{x}{y}_{1} &= \|y\|_{1}\sign(y)^\top x
      \end{aligned}$
      & $\begin{aligned}
	\mu_{1}(A) &= \max_{j \in \until{n}} \Big(a_{jj} + \sum_{i \neq j} |a_{ij}|\Big) \\ &= \sup_{\|x\|_{1} = 1}
	\sign(x)^\top Ax
      \end{aligned}$
      \\[2ex]
      &&  \\[-2ex]
      $ \begin{aligned}
	\norm{x}{\infty}
      \end{aligned}$
      &
      $\begin{aligned}
	\WP{x}{y}_{\infty} &= \max_{i \in I_{\infty}(y)} x_iy_i
      \end{aligned}$
      & $\begin{aligned}
	\mu_{\infty}(A) &= \max_{i \in \until{n}} \Big(a_{ii} + \sum_{j \neq i} |a_{ij}|\Big) \\ &=\max_{\|x\|_{\infty} = 1}
	\max_{i \in I_{\infty}(x)} (Ax)_ix_i
      \end{aligned}$
      \\[2ex] \hline
      && \\[-1ex]
  \end{tabular}}
  \caption{Weak pairings, and log norms for weighted $\ell_2$ and general
    $\ell_p$, for $1\leq p\leq\infty$ norms.  We adopt the shorthand
    $\Iinfty(x) = \setdef{i\in\until{n}}{|x_i|=\norm{x}{\infty}}$. The
    matrix $P$ is positive definite; the symbol $\circ$ stands for the
    Hadamard (element-wise) product. Only the unweighted $\ell_p$ norms,
    weak pairings, and log norms for $p \neq 2$ are included here since
    $\mu_{p,R}(A) = \mu_{p}(RAR^{-1})$ and $\WP{x}{y}_{p,R}=\WP{Rx}{Ry}_p$.
  }
    \label{table:equivalences}
\end{table*}


As will be discussed below, a special role is played by the \emph{sign
pairing} $\WP{\cdot}{\cdot}_{1}$ and \emph{max pairing}
$\WP{\cdot}{\cdot}_{\infty}$, associated respectively with the $\ell_1$ and
$\ell_{\infty}$ norms~\cite{AD-SJ-FB:20o}:
\begin{align}
  \label{def:WP-1+infty}
  \WP{x}{y}_{1} &= \|y\|_{1}\sign(y)^\top x , \\
  \WP{x}{y}_{\infty} &= \max_{i \in I_{\infty}(y)}  y_ix_i=\|y\|_{\infty}\max_{i \in I_{\infty}(y)}  \sign(y_i)x_i,
\end{align}
where $\Iinfty(y):= \setdef{i\in\until{n}}{|y_i|=\norm{y}{\infty}}$ is the
set of indices where $y$ takes its maximal absolute value.

While it is known from \cite{AD-SJ-FB:20o} that both the sign and max
pairings satisfy various properties of a regular pairing, as guaranteed by
Theorem~\ref{thm.main}, each statement in \cite{AD-SJ-FB:20o} is proved
using ad-hoc methods. The approach in this paper clearly demonstrates that
the sign and max pairings are regular pairings, and thus automatically
satisfy the numerous equivalent properties in Theorem~\ref{thm.main}. This
improves upon the original analysis in \cite{AD-SJ-FB:20o} and provides a
more streamlined and comprehensive treatment of the standard pairings
summarized in Table~\ref{table:equivalences}.

\begin{lemma}[\bf Regularity of Standard Pairings]\label{lem.regularity-standard}
  Each pairing collected in Table~\ref{table:equivalences} is regular. If
  $1<p<\infty$, then the norm $\norm{\cdot}{p,R}$ admits no regular pairing
  other than $\WP{x}{y}_{p,R}=\WP{Rx}{Ry}_p$ (being both LG and upper JMT
  pairings).
\end{lemma}
\begin{proof}
  It is known from~\cite{AD-SJ-FB:20o} that the binary operations
  $\WP{\cdot}{\cdot}_{p,R}$ are weak pairings.  A straightforward
  computation shows that all of them satisfy the straight angle property
  $\WP{-x}{x}=-\|x\|^2$. The final statement follows directly from
  Lemma~\ref{thm.unique}, noting that $\ell_p$ and weighted $\ell_p$ norms
  for $p\in(1,\infty)$ are $C^1$-smooth (and thus Gateaux differentiable)
  at every point except the origin.
\end{proof}

Since $\ell_1$ and $\ell_{\infty}$ norms are not Gateaux differentiable, Lemma~\ref{thm.unique} predicts that the regular pairing is non-unique.
Note that the sign-pairing $\WP{\cdot}{\cdot}_1$ is a special case of the LG pairing. It is also known from \cite[Example~13.1(b)]{KD:85} that the upper JMT pairing
(referred to as the ``Deimling pairing'' in~\cite{AD-SJ-FB:20o}) for the
$\ell_1$ norm is:
\begin{equation}\label{eq:ell1Deimling}
  [x,y]_{+,1} = \|y\|_{1} \Bigl(\sign(y)^\top x + \sum\nolimits_{i = 1}^n|x_i|\chi_{\{0\}}(y_i)\Bigr)\geq\WP{x}{y}_1,
\end{equation}
where $\chi_{0}(y_i)$ is the indicator function, equal to $1$ when $y_i=0$ and $0$ otherwise. Although, in accordance with Lemma~\ref{thm.unique},
$[\cdot,y]_{+,1} = \WP{\cdot}{y}_1$ for almost all $y$ (specifically, when $y\ne 0$), the sign pairing is more convenient for the curve norm derivation formula due to its simpler representation and its linearity in the first argument.

Similarly, the max pairing is not the unique regular pairing associated with the $\ell_{\infty}$ norm. The max pairing is nonlinear in both arguments, e.g.,
if $y=[1,\ldots,1]^{\top}$, then $\WP{x}{y}_{\infty}=\max_i x_i$. Hence, it is \emph{not} an LG pairing. However, one can easily construct LG pairings for $\ell_{\infty}$, norms; for instance, let
\[
\WP{x}{y}_{\infty,{\rm min-index}}=x_{m(y)}y_{m(y)}=\|y\|_{\infty}x_{m(y)}\sign y_{m(y)},
\]
where $m(y)=\min I_{\infty}(y)$ is the \emph{minimal} index $m$ such that $|y_m|=\|y\|_{\infty}$. Instead of the minimal index,
one can consider the maximal index. More generally, one can choose $m(y)=\mathcal{M}(I_{\infty}(y))$, where $\mathcal{M}:2^{\{1,\ldots,n\}}\to\{1,\ldots,n\}$
is a selector map, such that $\mathcal{M}(I)\in I$ for all non-empty $I\subseteq\{1,\ldots,n\}$, ensuring that $|y_{m(y)}|=\|y\|_{\infty}$.

An important advantage of the max pairing is that it is \emph{permutation-invariant} (remains unchanged under the permutation of coordinates).
In contrast, this property is clearly lost by the ``min-index'' LG pairing $\WP{x}{y}_{\infty,{\rm min-index}}$.

Unlike the sign pairing, the max pairing is a special case of the JMT pairing.
\begin{lemma}[\bf Max Pairing is the Upper JMT Pairing]\label{lem.max-pairing-jmt}
For any $x,y\in\mathbb{R}^n$,
\begin{equation}\label{eq.max-pairing-jmt}
\WP{x}{y}_{\infty}=[x,y]_{\infty,+}:=\|y\|_{\infty}\lim_{h\to 0+}\frac{\|y+hx\|_{\infty}-\|y\|_{\infty}}{h}.
\end{equation}
\end{lemma}
The proof is given in Appendix~\ref{app.proof-lem-max-pairing}. This proof trivially extends to the weighted max pairing;  we will consider a more general class of pairings in the next subsection.
	It should also be noticed that the infinite-dimensional counterpart of the relation~\eqref{eq.max-pairing-jmt}, dealing with continuous functions on a line segment $I=[a,b]$ and the standard norm $\|f\|_C=\max_{s\in I}|f(s)|$,  appears in~\cite[Theorem~4.5]{RCJ:47}. The latter theorem states that the upper JMT pairing, corresponding to this norm, is found as
	\[
	\begin{gathered}
		[g,f]_{C,+}=\|f\|_C\max_{s\in A(f)}g(s)\sign f(s)=\max_{s\in A(f)} g(s)f(s),\\
		A(f):=\{s\in I:|f(s)|=\|f\|_C\}.
	\end{gathered}
	\]

\subsection{Regular Pairings for Polyhedral Norms}\label{subsec:polyhedral-max}

The previously considered $\ell_1$ and $\ell_{\infty}$ norms are examples
of \emph{polyhedral norms}. In this subsection, we consider a more general
and novel case of a polyhedral norm, called the \emph{polyhedral max norm},
and study the induced norm and a \emph{polyhedral max pairing} related to
this norm.
\begin{definition}[Polyhedral and polyhedral max norms]
  A norm $\norm{\cdot}{}$ is \emph{polyhedral} if its unit disk
  $\setdef{v\in\real^n}{\norm{v}{}\leq 1}$ is a polyhedron. Given a full-column rank matrix $W\in\real^{m\times{n}}$ with $m\geq{n}$,
  the \emph{polyhedral max norm} $\norm{\cdot}{W}$ on $\real^n$ is defined
  by
  \begin{equation}
    \norm{x}{W} = \norm{Wx}{\infty} = \max_{i\in\until{m}} |w_i^\top x|,
  \end{equation}
  where $w_i\in\real^n$ is the $i$th row of $W$ (regarded as a column
  vector).
\end{definition}
Clearly, the standard non-Euclidean norms $\ell_1$ and $\ell_\infty$ are
polyhedral.  It is simple to verify that the polyhedral max norm
$\norm{\cdot}{W}$ satisfies the three defining properties of a norm (the
positive definite property follows from $W$ being full column
rank). Finally, from the hyperplane characterization of convex polytopes,
it is also easy to see that each polyhedral norm can be written as a
polyhedral max norm.

\begin{definition}
  The \emph{polyhedral max pairing}
  ${\WP{\cdot}{\cdot}_W}$ 
  is defined by
  \begin{equation}
    \WP{x}{y}_W =  \WP{Wx}{Wy}_\infty\quad\forall x,y\in\real^n.
  \end{equation}
\end{definition}

We are now ready to state our main result on polyhedral pairings.

\begin{lemma}[\bf Properties of Polyhedral Max Pairings] \label{lem:polyh-max-pairing}
  For a full-column rank matrix $W\in\real^{m\times{n}}$ with $m\geq{n}$,
  the polyhedral max pairing $\WP{\cdot}{\cdot}_W$ is the upper JMT pairing
  (and hence a regular pairing), associated with the norm
  $\norm{\cdot}{W}$.

  For each $A\in\real^{n\times{n}}$, its associate log norm is
  \begin{subequations} \label{opt-problem-lognorm}
    \begin{align}
      \mu_W(A) \quad = \quad  \min_{H\in\real^{m\times{m}}} \quad & \mu_\infty(H) \\
      \subject \quad
      &WA  = H W, \label{eq:weird-comm}
    \end{align}
  \end{subequations}
  where the linear constraint~\eqref{eq:weird-comm} is feasible and the
  minimum exists\footnote{Notice that the set of admissible matrices $H$ is
  an affine subspace in $\mathbb{R}^{m\times m}$, which is non-compact
  (closed yet unbounded), and the convex function $\mu_{\infty}$ is not
  radially unbounded. Hence, the existence of a minimum
  in~\eqref{opt-problem-lognorm} is a self-standing result that is not
  fully obvious.}.
\end{lemma}
The proof of Lemma~\ref{lem:polyh-max-pairing} will be given in Appendix~\ref{app.proof-lem-polyh-max-pairing}.

It should be noted that the convex program~\eqref{eq:weird-comm} can, in
fact, reduces\footnote{It should be emphasized, however, that matrix $W$
needs to be fixed. As discussed by~\cite{FB:99}, when only $A$ is given and
one is interested to find the polyhedral norm with minimal value of
$\mu_W(A)$, the constraint~\eqref{eq:weird-comm} is bilinear in the
unknowns $(W,\barAA)$ so the feasible set is non-convex.} be written as a
standard LP with $2m^2-m$ variables, as implied by the following remark.

\begin{remark} When $m>n$, the matrix $\barAA\in\real^{m\times{m}}$ is
under-determined by the $m\times{n}$ equalities $WA = \barAA W$, in other
words, the solutions to equation~\eqref{eq:weird-comm} live in an
$m(m-n)$-dimensional vector subspace. Visually,
\begin{equation*}
  \renewcommand\matscale{.45}
  \matbox{8}{4}{m}{n}{W}
  \raiserows{2.0}{\matbox{4}{4}{n}{n}{A}}
  \enspace=\enspace \matbox{8}{8}{m}{m}{H}
  \raiserows{0}{\matbox{8}{4}{m}{n}{W}}.  
\end{equation*}
\end{remark}
Also note that, at fixed $A$ and $W$, the optimization
problem~\eqref{opt-problem-lognorm} can be rewritten as a linear program by
parametrizing the variable $H$ as follows.  From~\cite{SJ-AD-AVP-FB:21f},
we know that $\lognorm{H}{\infty}\leq\gamma$ if and only if there exists
zero-diagonal $T\in\real^{n\times{n}}$ satisfying
\begin{subequations}\label{eq:weird-comm-equiv}
  \begin{gather}
    \diag(H) - T \leq H \leq \diag(H) + T,  \label{fact:polytope-lognorm:2:1}\\
    (\diag(H)+T) \vect{1}_n \leq \gamma \vect{1}_n. \label{fact:polytope-lognorm:2:2}
  \end{gather}
\end{subequations}
where $\diag(H)$ is the diagonal matrix with entries equal to the diagonal
entries of $H$.  In this transcription we used the fact that, for a Metzler
$M$, $\mu_\infty(M)=\min\setdef{b\in\real}{M\vect{1}_n\leq b \vect{1}_n}$;
e.g., see~\cite[Lemma~2.8]{FB:24-CTDS}.  Hence, finding the minimum
in~\eqref{opt-problem-lognorm} is equivalent to finding the minimal
$\gamma$, for which constraints~\eqref{eq:weird-comm-equiv} are feasible.
\oprocend

It can be noticed that the $\ell_1$ norm is a special case of polyhedral
max norm, where matrix $W$ has $m=2^n$ rows, among which are all possible
$n$-tuples of $\pm 1$. Lemma~\ref{lem:polyh-max-pairing} implies that the
corresponding polyhedral max pairing coincides
with~\eqref{eq:ell1Deimling}, which can also be verified directly through
some tedious computations.

\section{Contraction Analysis Based upon Regular Pairings}\label{subsec:motivating-section}

In this section, we review the basic features of contraction theory in light of our results on regular pairings, specifically the curve norm derivative formula and Lumer's property.

We begin by illustrating some direct implications of the curve norm derivative formula and Lumer's property.
Given a curve
  $\map{z}{\real_{\geq0}}{\real^n}$ and a norm with a compatible regular pairing
  $\WP{\cdot}{\cdot}$, the curve norm derivative formula implies
  \begin{equation}
    \tfrac{1}{2}\tfrac{d}{dt}\norm{z(t)}{}^2
    = \WP{\dot{z}(t)}{z(t)}.
  \end{equation}
  for almost all $t$. For instance, given a trajectory $z(t)$ that is a solution to a time-varying system $\dot z(t)=A(t)z(t)$ with a locally bounded matrix $A(\cdot)$, the curve norm derivative formula and Lumer's property imply
\begin{align}
\tfrac{1}{2}\tfrac{d}{dt} \norm{z(t)}{}^2 = \WP{A(t)z(t)}{z(t)} \leq \lognorm{A(t)}{} \norm{z(t)}{}^2.
\end{align}
Hence, upper bounds on the matrix's induced norm can be easily transformed into explicit estimates of the solution's norm.
For commonly used norms, we have:
\begin{alignat*}{3}
   \lognorm{A(t)}{1}&\leq b \quad&&\implies\quad \tfrac{1}{2} \tfrac{d}{dt} \norm{z}{1}^2 &=&\,
   \norm{z}{1} \sign(z)^{\top} \dot{z} \leq b \norm{z}{1}^2, \\
   \lognorm{A(t)}{2}&\leq b \quad&&\implies\quad
   \tfrac{1}{2} \tfrac{d}{dt}  \norm{z}{2}^2 &=&\,  z^\top \dot{z}
    \leq b \norm{z}{2}^2,
   \\
   \lognorm{A(t)}{\infty}&\leq b \quad&&\implies\quad
   \tfrac{1}{2} \tfrac{d}{dt}  \norm{z}{\infty}^2 &=&\,
   \max_{i\in \Iinfty(z)} \left\{z_i \dot{z}_i  \right\}
   \leq b \norm{z}{\infty}^2.
 \end{alignat*}
 The importance of these three norms (and their weighted counterparts) lies in the possibility of explicitly computing the log norm.
 The same applies to the polyhedral max norm, examined in the previous section. Unlike these, no algorithm for computing $\mu_p(A)$ for $p\ne 1,2,\infty$ is known to the best of the authors' knowledge. Remarkably, the computation of the induced operator $\ell_p$-norm is known to be NP-hard~\cite{AB-AV:11}.

\subsection{Contraction Criterion}

Contraction analysis deals with a similar argument, where $z(t)=x_1(t)-x_2(t)$ is a discrepancy of two solutions to the same system
\[
\dot x(t)=f(t,x(t)).
\]
Assuming the right-hand side to be $C^1$-smooth in $x$, $z(t)$ satisfies the equation
\[
\dot z(t)=A(t)z(t),\quad A(t)=\int_0^1\frac{\partial f}{\partial x}(\theta x(t)+(1-\theta)y(t))d\theta.
\]

Although the matrices $A(t)$ usually cannot be computed explicitly, they belong to the convex hull spanned by the Jacobian values, which often allows estimation of their log norm $\mu(A(t))$ using the integral Jensen inequality. Our characterization of regular pairings (Theorem~\ref{thm.main}) leads to the following result, which clarifies the contraction criterion for continuously differentiable vector fields\footnote{A similar result can also be proved under local Lipschitz continuity in $x$; see~\cite{AD-AVP-FB:22q}.} as presented in \cite{AD-SJ-FB:20o}.

\begin{theorem}[\bf Contraction Criteria for $C^1$-Smooth Vector Fields]
  \label{thm:general}
  Consider the dynamics $\dot{x} = f(t,x)$, with $f$ continuously
  differentiable in $x$ and continuous in $t$. Let $C \subseteq \real^n$ be
  an open, convex, and forward invariant set and let $\|\cdot\|$ denote a norm
  with a compatible regular pairing $\WP{\cdot}{\cdot}$. Then, for $b \in
  \real$, the following statements are equivalent:
  \begin{enumerate}
  \item\label{ctGen:5} $\WP{f(t,x) - f(t,y)}{x - y} \leq b\|x - y\|^2$, for all $x,y \in C, t \geq 0$;
  \item\label{ctGen:4} $\WP{\jac{f}(t,x)v}{v} \leq b\|v\|^2$, for all $v \in \real^n, x \in C, t \geq 0$;
  \item\label{ctGen:2} $\mu(\jac{f}(t,x)) \leq b$, for all $x \in C, t \geq 0$;
  \item\label{ctGen:6} $D^+\|\phi(t,t_0,x_0) - \phi(t,t_0,y_0)\| \leq b\|\phi(t,t_0,x_0) - \phi(t,t_0,y_0)\|$, for all $x_0,y_0 \in C, 0 \leq t_0 \leq t$ for which the solutions exist\footnote{Here $D^+$ stands for the right upper Dini derivative, see Eq.~\eqref{eq.dini-def} in Appendix.};
  \item\label{ctGen:1} $\|\phi(t,t_0,x_0) - \phi(t,t_0,y_0)\| \leq e^{b(t-s)}\|\phi(s,t_0,x_0) - \phi(s,t_0,y_0)\|$, for all $x_0, y_0 \in C$ and $0 \leq t_0 \leq s \leq t$ for which the solutions exist.
  \end{enumerate}
  Here $\jac{f}(t,x) := \frac{\partial f}{\partial x}(t,x)$ is the Jacobian
  of $f$.
\end{theorem}
\begin{proof}
According to~\cite[Theorem~31]{AD-SJ-FB:20o}, the conditions~{\it\ref{ctGen:5}-\ref{ctGen:1}} are equivalent for any weak pairing that satisfies the curve norm derivative formula and the rightmost inequality in~\eqref{eq.dominated} (referred in~\cite{AD-SJ-FB:20o} to as the ``Deimling inequality''). According to our Theorem~\ref{thm.main}, this holds if and only
if the pairing is regular.
\end{proof}
\begin{remark}
Referring to Table~\ref{table:equivalences}, one observes that when the norm is the Euclidean weighted norm induced by a positive definite matrix $P \succ 0$, i.e., $\|x\|_{2,P^{1/2}} = (x^\top P x)^{1/2}$, condition (iii) reduces to the well-known Demidovich condition~\cite{AP-AP-NVDW-HN:04}:
\begin{equation}\label{eq:demid}
P\jac{f}(t,x)+\jac{f}(t,x)^{\top}P\preceq 2bP\quad\forall t\geq 0,x\in C.
\end{equation}
The seminal work by~\cite{WL-JJES:98} extends this criterion to the case of Riemannian norms, where $P=P(x)$ is a Riemannian matrix defining a pointwise inner product on the tangent space. Theorem~\ref{thm:general}, by contrast, extends the Demidovich criterion to non-Euclidean norms—most notably, to weighted $\ell_1$ and $\ell_\infty$ norms.
\end{remark}
\begin{remark}
  The condition~{\it\ref{ctGen:5}} is known as the \emph{one-sided Lipschitz condition}, which is most typically validated with respect to
  an inner product or the JMT upper pairing~\cite{GS:24}. Theorem~\ref{thm:general} shows that, to prove the contraction property,
  the JMT upper pairing can be replaced by an arbitrary regular pairing.\oprocend
\end{remark}

\subsection{Example}
We consider the prototypical biochemical reaction model studied by~\cite{GR-MDB-EDS:10a,ZA-EDS:13,ZA-EDS:14b}, in which an enzyme $X$ binds with a substrate $S$ to form a complex $Y$~\cite{ZA-EDS:13}. Let $x(t)$ and $y(t)$ denote the concentration of the enzyme and the complex, respectively. The dynamics are given by
\begin{equation}\label{eq:biochem}
  \begin{bmatrix} \dot{x} \\  \dot{y} \end{bmatrix}
  = \begin{bmatrix}
    -k_0 x + k_1 y -k_2 (\ymax-y) x + u \\
    - k_1 y + k_2 (\ymax-y) x
  \end{bmatrix} ~:= f(x,y) + \begin{bmatrix} u \\ 0
  \end{bmatrix},
\end{equation}
where $k_0>0$ is a degradation rate, $k_1,k_2>0$ are reaction rates,
$\ymax>0$ is a concentration of the substrate (considered as constant), and $u=u(t)$ is an external
input. It can be easily shown~\cite{FB:24-CTDS} that the convex set $C=\{(x,y):x\geq 0,0\leq y\leq\ymax\}$ is forward invariant for every function $u(t)\geq 0$.

\textbf{Contraction in a polyhedral norm.}
As demonstrated in~\cite[Section~4.6]{FB:24-CTDS},
the system~\eqref{eq:biochem} proves to be contractive in the weighted $\ell_1$-norm $\|\cdot\|_{1,R}$, where
\[
R(q)=\begin{bmatrix}
1 & 0\\
0 & q
\end{bmatrix},\quad\text{with}\quad 1<q<1+\frac{k_0}{k_2\ymax}.
\]
Indeed, computing the Jacobian $\jac{f}(x,y)$ and noticing that
   \begin{equation}
     R(q)\jac{f}(x,y)R(q)^{-1} = \begin{bmatrix}
       -k_0 - k_2 (\ymax-y) & (k_1+k_2 x)q^{-1} \\
       qk_2 (\ymax-y) & -k_1-k_2 x
     \end{bmatrix},
   \end{equation}
one easily proves that, for each $x\geq 0$ and $0\leq y\leq\ymax$,
\[
\begin{aligned}
\mu_{1,R(q)}(\jac{f}(x,y))&=\max\left\{(1-q^{-1})(k_1+k_2x),-k_0+(q-1)(\ymax-y)\right\}\\
&\leq b(q):=\max\{k_1-k_1q^{-1},-k_0+(q-1)\ymax\}<0.
\end{aligned}
\]

It is noteworthy that, as demonstrated in \cite{ZA-EDS:13}, the system fails to be contractive in any diagonally weighted $\ell_p$ norm with $p>1$: for every $p>1$ and every positive diagonal matrix $R$, there exists some $(x,y)\in C$ such that
$\mu_p\bigl(R\,D f(x,y)\,R^{-1}\bigr) \;>\; 0$.

\textbf{Nonexistence of Quadratic Lyapunov Functions.} We will now show
that for sufficiently large $\ymax$, the system cannot be contractive -- or
even non-expansive -- in any Euclidean norm. Specifically, we will show
that inequality \eqref{eq:demid} fails at some point $(x,y)\in C$ for every
$P\succ0$, even when $b=0$. Indeed, introducing the new variable
$z:=\ymax-y\in [0,\ymax]$, one can compute
$P\jac{f}(x,y)+\jac{f}(x,y)^{\top}P$ to be
\[
\begin{bmatrix}
-2k_0 p_{11} - 2k_2z p_{11} + 2k_2z p_{12}
&
(p_{11} - p_{12})(k_1 + k_2 x) - k_0 p_{12} - k_2zp_{12} + k_2zp_{22}
\\[1.2ex]
*
&
2(p_{12} - p_{22})(k_1 + k_2 x)
\end{bmatrix}.
\]
One first notices that, for this matrix to be negative semidefinite as $x\to+\infty$, the leading term that depends linearly on $x$ must be negative semidefinite
\[
\begin{bmatrix}
0 & (p_{11}-p_{12})k_2x\\
(p_{11}-p_{12})k_2x & 2(p_{12}-p_{22})x
\end{bmatrix}
\preceq0
\]
for large~$x$. By Sylvester's criterion, this can occur only if $p_{11}=p_{12}$ and $p_{12}\leq p_{22}$.
Since $P\succ 0$, the latter inequality must be strict: $p_{12}<p_{22}$.
Assuming that these assumptions hold, one can further simplify the matrix in~\eqref{eq:demid} as
\[
P\jac{f}(x,y)+\jac{f}(x,y)^{\top}P=M(x,z):=\begin{bmatrix}
-2k_0 p_{11}
&
- k_0 p_{11} + k_2z(p_{22}-p_{11})
\\[1.2ex]
*
&
2(p_{11} - p_{22})(k_1 + k_2 x)
\end{bmatrix}.
\]
The determinant of this matrix can be written as a quadratic formula:
\[
\begin{gathered}
\det M(x,z)=\alpha z^2+\beta z+\gamma(x),\\
\alpha := -k_2^2\,(p_{22}-p_{11})^2<0,\;\;\beta := 2\,k_0\,k_2\,p_{11}\,(p_{22}-p_{11})>0,\\
\gamma(x) := 4\,k_0\,p_{11}\,(p_{22}-p_{11})\,(k_1 + k_2 x)\;-\;k_0^2\,p_{11}^2.
\end{gathered}
\]
We claim that the latter expression can be nonnegative for $x=0$ and all
$z\in [0,\ymax]$ only when $\ymax$ is small enough.  Indeed, first,
$M(0,0)\geq 0$ if and only if
\[
\gamma(0)=4\,k_0\,p_{11}\,(p_{22}-p_{11})\,k_1-k_0^2\,p_{11}^2\geq 0\Longrightarrow
\frac{p_{11}}{p_{22}-p_{11}}\leq \frac{4k_1}{k_0}.
\]
Second, $\ymax$ should not exceed the maximal of the two roots of $M(0,z)$, i.e.,
\[
\ymax\leq\frac{-\beta-\sqrt{\beta^2-4\alpha\gamma}}{2\alpha}=\frac{p_{11}}{p_{22}-p_{11}}\frac{k_0}{k_2}
+2\sqrt{\frac{p_{11}}{p_{22}-p_{11}}\frac{k_0k_1}{k_2^2}}=\frac{8k_1}{k_2}.
\]
Hence, for $\ymax>8k_1/k_2$ (the concentration of the substrate is large
enough), the system is not contractive in any norm $\|\cdot\|_{2,P^{1/2}}$.

\section{Conclusion}\label{sec:concl}

This paper introduces a notion of \emph{regular pairing}, along with a characterization theorem that demonstrates the equivalence of several useful properties
and provides a computationally-friendly regularity criterion for weak pairings. We describe how the concept of a regular pairing is weaker than previous notions, retaining the critical properties of Lumer's equality and the curve norm derivative formula. These tools are important for Lyapunov stability analysis and contraction theory, allowing for elegant and compelling proofs. Since regular pairings are unique for smooth norms, we give special attention to non-differentiable polyhedral norms, defining and characterizing a polyhedral max pairing for the polyhedral max norm.
In conclusion, we outline several directions for future research.

Non-Euclidean norms -- such as \emph{polynomial}~\cite{AAA-EdK-GH:19} and \emph{Barabanov} norms~\cite{RT-MM:12} -- often arise in Lyapunov-based stability analysis of switching systems. Analyzing regular pairings associated with these norms, in light of Theorem~\ref{thm:general}, is expected to yield new stability and contraction criteria for switching systems with external inputs, as well as for interconnections of switching systems~\cite{XD-SJ-FB:19f,FB:24-CTDS}.
The potential of \emph{polyhedral} Lyapunov functions -- demonstrated, for example, in recent work on biochemical reaction networks~\cite{MAAR-DA-EDS:23} -- remains underexplored, largely due to the absence of efficient design methods~\cite{MAAR-DA-EDS:23}. Given a polyhedral max norm, the associated log norm -- which plays a central role in establishing contraction -- can be computed via Lemma~\ref{lem:polyh-max-pairing} by solving the convex optimization problem~\eqref{opt-problem-lognorm}. A challenging direction for future research is the design of such norms, specifically the construction of a matrix $W$ of minimal dimension that ensures $\mu_W(A) \leq c < 0$ for a given family of matrices $A$. Finally, Lumer's work on semi-inner products and the numerical ranges of linear operators was, to a great extent, motivated by the analysis of stability in linear time-invariant systems within Banach spaces (more specifically, by the study of contraction properties of $C_0$-semigroups~\cite{GL-RSP:61}). As implied by Corollary~\ref{cor.lg-sip}, LG pairings form a special subclass of regular pairings in a real normed space; however, their construction for a general norm is not explicit, as it relies on the Hahn-Banach theorem. An important direction for ongoing research is to identify explicit regular pairings for functional norms and to develop corresponding contraction criteria -- along the lines of Theorem~\ref{thm:general} -- for nonlinear partial differential equations and integral equations.

\section*{Acknowledgments}

The authors thank Alexander Davydov for insightful conversations on pairings and contraction theory.

\bibliographystyle{siamplain}
\bibliography{alias,Main,FB}
\clearpage
\appendix
\section{A Lemma on the Differentiation of Curve Norms}\label{app.proof-lem-curvenorm}
In this section, we prove a technical lemma that establishes a strengthened version of the curve norm derivative formula for the upper JMT pairing, generalizing the result of~\cite[Proposition~13.1]{KD:85} to arbitrary normed spaces (which may be incomplete).

Given a function $f:(a,b)\to\real$, one can define the quadruple of \emph{Dini derivatives} at $t_0\in (a,b)$ as follows
\begin{equation}\label{eq.dini-def}
\begin{gathered}
D_+f(t_0)=\liminf_{h\to 0+}\frac{f(t_0+h)-f(t_0)}{h},\;\;D^+f(t_0)=\limsup_{h\to 0+}\frac{f(t_0+h)-f(t_0)}{h},\\
D_-f(t_0)=\liminf_{h\to 0-}\frac{f(t_0+h)-f(t_0)}{h},\;\;D^-f(t_0)=\limsup_{h\to 0-}\frac{f(t_0+h)-f(t_0)}{h}.
\end{gathered}
\end{equation}
Obviously, if the derivative $\dot f(t_0)$ exists, it coincides with all the Dini derivatives.

\begin{lemma}[\bf Dini Derivatives and the JMT Pairings]\label{lem.curve-norm-deriv}
  For every differentiable curve $x:(a,b)\to X$ and $t\in (a,b)$ the
  following statements are valid:
  \begin{enumerate}
  \item all Dini derivative of the function $f(t):=\|x(t)\|$ are finite, and
    \begin{equation}\label{eq.dini-bound}
      -\|\dot x(t)\|\leq D_-f(t),D_+f(t),D^-f(t),D^+f(t)\leq \|\dot x(t)\|;
    \end{equation}
  \item the Dini derivatives are related to the JMT pairings by the equalities
    \begin{equation}\label{eq.deimling}
      \begin{gathered}
        f(t)\,D^+f(t)=f(t)\,D_+f(t)=\SIP{\dot x(t)}{x(t)}_+\qquad\forall t\in (a,b),\\
        f(t)\,D^-f(t)=f(t)\,D_-f(t)=\SIP{\dot x(t)}{x(t)}_-\qquad\forall t\in (a,b).
      \end{gathered}
    \end{equation}
    In particular, if $f(t)\ne 0$, then the right and left derivatives $\dot
    f(t\pm 0)$ exist;\\
  \item the derivative $\dot f(t)$ exists at almost all $t\in(a,b)$.
  \end{enumerate}
\end{lemma}
\begin{proof}

To prove statement~\emph{(i)}, notice that $|f(t+h)-f(t)|=|\,\|x(t+h)\|-\|x(t)\|\,|\leq \|x(t+h)-x(t)\|$ for each $h$. Hence,
\[
\frac{|f(t+h)-f(t)|}{|h|}\leq\left\|\frac{x(t+h)-x(t)}{h}\right\|\qquad\forall h\ne 0.
\]
Taking $\liminf$ and $\limsup$  as $h\to 0+$ and as $h\to 0-$, one proves~\eqref{eq.dini-bound}.

The last statement~\emph{(iii)} now follows from the celebrated
Denjoy-Young-Saks theorem~\cite[Chapter~IV, Theorem~4.4]{AMB:78} stating
that at almost every point $t\in(a,b)$ either $\dot f(t)$ exists or at
least two of the four Dini derivatives are infinite (which, as we have
shown, is impossible in our case).

The statement~\emph{(ii)} is trivial in the case where $x(t)\ne\vect{0}$ (i.e., $f(t)=0$). Indeed, the Dini derivatives are all finite~\eqref{eq.dini-bound}, and hence $f(t) D^{\pm}f(t)=f(t) D_{\pm}f(t)=0=\SIP{\dot x}{\vect{0}}$.
If $x(t)\ne\vect{0}$, recall that $x(t+h)=x(t)+h\dot x(t)+v(t,h)$, where $h^{-1}\|v(t,h)\|\xrightarrow[h\to 0]{} 0$, by definition of the derivative. Therefore,
\[
f(t)\frac{f(t+h)-f(t)}{h}=\|x(t)\|\frac{\|x(t)+h\dot x(t)\|-\|x(t)\|}{h}+\alpha(t,h),\;\;
\]
where $\alpha(t,h)$ vanishes as $h\to 0$. Passing to the limit as $h\to 0+$ and using~\eqref{eq.JMT-SIP} (where $x,y$ are replaced by, respectively, $\dot x$ and $x$), one shows that
\begin{equation}\label{eq.aux2}
\lim_{h\to 0+}f(t)\frac{f(t+h)-f(t)}{h}=\SIP{\dot x(t)}{x(t)}_+\qquad\forall t\in (a,b).
\end{equation}
This, obviously, entails the first line in equations~\eqref{eq.deimling}:
\[
\dot f(t+0)=D^+f(t)=D_+f(t)=\lim_{h\to 0+}\frac{f(t+h)-f(t)}{h}\overset{\eqref{eq.aux2}}{=}\frac{\SIP{\dot x(t)}{x(t)}_+}{f(t)},
\]
The second line in~\eqref{eq.deimling} is proved in the same way by considering limits as $h\to 0-$.
\end{proof}

\section{Proof of Theorem~\ref{thm.main}}\label{app.proof-thm}

We first prove the equivalence of statements~\emph{\ref{thm1:straight}-\ref{thm1:lumer1}}.

\textbf{Implication $\ref{thm1:straight}\Longrightarrow\ref{thm1:linear}$.} Combining the straight angle property with~\eqref{eq.weak-hom1}, one shows that $\WP{ay}{y}=a\|y\|^2$ for all $a\in\real$. For $a\geq 0$ this is implied by~\eqref{eq.weak-hom1}. If $a=-|a|<0$, then
$\WP{ay}{y}=|a|\WP{-y}{y}=-|a|\|y\|^2=a\|y\|^2$. Thanks to~\eqref{eq.subadd},
\[
\WP{x}{y}=\WP{x+ay-ay}{y}\overset{\eqref{eq.subadd}}{\leq} \WP{x+ay}{y}-a\|y\|^2\overset{\eqref{eq.subadd}}{\leq}
\WP{x}{y}+a\|y\|^2-a\|y\|^2=\WP{x}{y}
\]
for all $x,y\in X$, which means that $\WP{x}{y}+a\|y\|^2=\WP{x+ay}{y}$.

\textbf{Implication $\ref{thm1:linear}\Longrightarrow\ref{thm1:dominated}$.} Assuming that the WP is partially linear and fixing $x,y\in X$, $h>0$ one obtains, due to the Cauchy-Schwarz inequality, that
\[
\WP{x}{y}=\frac{\WP{y+hx}{y}-\|y\|^2}{h}\overset{\eqref{eq.cauchy-wp}}{\leq}\frac{\|y+hx\|\,\|y\|-\|y\|^2}{h}=\|y\|\frac{\|y+hx\|-\|y\|}{h}.
\]
Passing to the limit as $h\to 0+$, one proves the rightmost inequality in~\eqref{eq.dominated}. To prove the remaining inequalities in~\eqref{eq.dominated}, it now suffices to notice that
\[
\SIP{x}{y}_-=-\SIP{-x}{y}_+\leq -\WP{-x}{y}\overset{(!)}{\leq}\WP{x}{y},
\]
where the inequality (!) is straightforward by substituting $x_1=x=-x_2$ into~\eqref{eq.subadd}.

\textbf{Implication $\ref{thm1:dominated}\Longrightarrow\ref{thm1:lumer1}$} is proved similarly to the first part of~\cite[Theorem~18]{AD-SJ-FB:20o}.
Notice that $\|x+hAx\|\leq\|\Id+hA\|$ whenever $A$ is bounded and $x$ is a unit vector ($\|x\|=1$).
Using the definition of the upper JMT pairing~\eqref{eq.JMT-SIP}, one obtains
\[
\WP{Ax}{x}\overset{\ref{thm1:dominated}}{\leq}\SIP{Ax}{x}_+=\lim_{h\to 0}\frac{\|x+hAx\|-\|x\|}{h}\leq \lim_{h\to 0}\frac{\|\Id+hA\|-1}{h}=\mu(A).
\]
whenever $\|x\|=1$, which is nothing else than the condition~\emph{\ref{thm1:lumer1}}.

\textbf{Implication $\ref{thm1:lumer1}\Longrightarrow\ref{thm1:straight}$} is proved by substituting $A=-\Id_X$ into~\eqref{eq.lumer1}:
\[
\WP{-x}{x}\leq -1\qquad\forall x\in X:\|x\|=1\overset{\eqref{eq.weak-hom1}}{\Longrightarrow} \WP{-x}{x}\leq-\|x\|^2\;\;\forall x\in X.
\]
In view of the Cauchy-Schwarz inequality, one also has $|\WP{-x}{x}|\leq\|x\|^2$, and thus the $\WP{-x}{x}=-\|x\|^2$, which finishes the first part of the proof.
\vskip0.3cm
Next, we show that the conditions~\emph{\ref{thm1:curve1},~\ref{thm1:curve2}} are equivalent to each other and to statements~\emph{\ref{thm1:straight}-\ref{thm1:lumer1}}. This is proved by the following three implications.

\textbf{Implication $\ref{thm1:curve1}\Longrightarrow\ref{thm1:curve2}$} is straightforward.

\textbf{Implication $\ref{thm1:curve2}\Longrightarrow\ref{thm1:straight}$.} Considering the linear curve $x(t)=x_0-tx_0=(1-t)x_0$ whose norm is differentiable at point $t=0$ and using~\eqref{eq.curve-norm},
one shows that
\[
\WP{-x_0}{x_0}=\WP{\dot x(0)}{x(0)}=\frac{1}{2}\biggl.\frac{d}{dt}\left(\|x_0\|^2(1-t)^2\right)\biggr|_{t=0}=-\|x_0\|^2\qquad\forall x_0\in X.
\]

\textbf{Implication $\ref{thm1:dominated}\Longrightarrow\ref{thm1:curve1}$.}
Suppose that \emph{\ref{thm1:dominated}} holds and consider a differentiable curve $x:(a,b)\to X$ whose norm $f(t)=\|x(t)\|$
 is differentiable at some $t_0\in (a,b)$. We are going to show that~\eqref{eq.curve-norm} holds at $t=t_0$. Without loss of generality, assume that $t_0=0$ and $b=-a>0$. In accordance with Lemma~\ref{lem.curve-norm-deriv},
\begin{equation}\label{eq.aux}
  f(0)\dot f(0)=f(0)D^+f(0)\overset{\eqref{eq.deimling}}{\geq}\SIP{\dot x(0)}{x(0)}_+\overset{\ref{thm1:dominated}}{\geq} \WP{\dot x(0)}{x(0)}.
\end{equation}
Consider now the function $\tilde x(t)=x(-t)$ whose norm $\tilde f(t)=\|\tilde x(t)\|=f(-t)$ is, obviously, also differentiable at $t=0$. Applying~\eqref{eq.aux} to $\tilde x,\tilde f$, one shows that
\begin{equation}\label{eq.aux+}
  -f(0)\dot f(0)=\tilde f(0)\dot{\tilde f}(0)\geq\WP{\dot{\tilde x}(0)}{\tilde x(0)}=\WP{-\dot x(0)}{x(0)}.
\end{equation}
Summing up the inequalities~\eqref{eq.aux} and~\eqref{eq.aux+} and using the subadditivity property~\eqref{eq.subadd},
\[
0=f(0)\dot f(0)-f(0)\dot f(0)\geq \WP{\dot x(0)}{x(0)}+\WP{-\dot x(0)}{x(0)}\overset{~\eqref{eq.subadd}}{\geq}
\WP{\vect{0}}{x(0)}=0.
\]
The latter inequality can only hold when $f(0)\dot f(0)=\WP{\dot x(0)}{x(0)}$. Hence,~\eqref{eq.curve-norm} holds at $t=t_0$ whenever $\dot f(t_0)$ exists, which finishes the proof of implication $\ref{thm1:dominated}\Longrightarrow\ref{thm1:curve1}$.
\vskip0.3cm
\textbf{The equivalence of~\emph{\ref{thm1:lumer}} and~\emph{\ref{thm1:straight}-\ref{thm1:curve2}}.} It remains to prove that the Lumer property is equivalent to the remaining conditions. The implication $\ref{thm1:lumer}\Longrightarrow \ref{thm1:lumer1}$ is straightforward, and hence~\emph{\ref{thm1:lumer}} implies all conditions~\emph{\ref{thm1:straight}-\ref{thm1:curve2}}.
To prove the opposite implication, it remains to prove that~\emph{\ref{thm1:straight}-\ref{thm1:curve2}} imply the inequality 
\begin{equation}\label{eq.aux1}
  \mu_0:=\sup_{x:\|x\|=1}\WP{Ax}{x}=\sup_{x\ne\vect{0}}\frac{\WP{Ax}{x}}{\|x\|}\geq \mu(A)\quad\forall A\in\mathcal{B}(X).
\end{equation}
To prove this, choose a vector $x_0\in X$ such that $\|x_0\|=1$ and consider the functions $x(t)=x_0+tAx_0$, $f(t)=\|x(t)\|$.
Applying~\eqref{eq.curve-norm1} to $v=Ax$, for almost all $t$ one has
\[
\begin{aligned}
f(t)\dot f(t)=\WP{Ax_0}{x(t)}=\WP{Ax(t)-tA^2x_0}{x(t)}\leq \WP{Ax(t)}{x(t)}+\WP{-tA^2x_0}{x(t)}\\
\overset{\eqref{eq.aux1},\eqref{eq.cauchy-wp},\,\|A^2x_0\|\leq c}{\leq} \mu_0\|x(t)\|+c|t|\,\|x(t)\|=(\mu_0+|t|c)f(t).
\end{aligned}
\]
Here $c:=\|A^2\|$ is a constant. Choosing now $h>0$ so small that $f(t)>0$ for $t\in [0,h)$, one has
$
\dot f(t)\leq \mu_0+tc$ for all $t\in [0,h).$
Notice that $f$ is Lipschitz, and hence also absolutely continuous. Integrating the inequality over $[0,h]$, one proves that
\[
\|x_0+hAx_0\|-1=f(h)-f(0)\leq\int_0^h(\mu_0+tc)dt=\mu_0h+\frac{h^2c}{2}\qquad\forall x_0:\|x_0\|=1.
\]
Taking the supremum  over all $x_0$ such that $\|x_0\|=1$, one shows that
\[
\|\Id_X+hA\|-1\leq\mu_0h+O(h^2),
\]
and therefore $\mu(A)=\lim_{h\to 0+}h^{-1}(\|\Id_X+hA\|-1)\leq\mu_0$, which, in view of~\emph{\ref{thm1:lumer1}}, proves the Lumer equality~\eqref{eq.lumer}. The proof of Theorem~\ref{thm.main} is finished.
\section{Proof of Lemma~\ref{thm.unique}}\label{app.proof-thm-unique}

The ``if'' part of the first statement is straightforward from Theorem~\ref{thm.main} and the result by Giles~\cite[Theorem~1]{JRG:67}, stating that every normed space admits at least one LG pairing. Indeed, if the lower and upper JMT pairings are equal, then every regular pairing is coincident with them in view of~\eqref{eq.dominated}. Taking the LG pairing $\SIP{\cdot}{\cdot}$ compatible with the norm, Corollary~\ref{cor.lg-sip} entails that $\realpart\SIP{x}{y}=\SIP{x}{y}_+$ for all $x,y\in X$. Using Lemma~\ref{lem.wp-and-sip}, one proves the uniqueness of such an LG pairing.

To prove the ``only if'' part, recall that $\SIP{x}{y}_-=-\SIP{-x}{y}_+$. Hence, if the normed space is not Gateaux differentiable,
one has $\SIP{x}{y}_+\ne-\SIP{-x}{y}_+$ for some pair of vectors $x,y\in X$. Obviously, such a regular pairing does not satisfy the conditions of Lemma~\ref{lem.wp-and-sip}, and thus cannot be represented as a real part of an LG pairing. According to the Giles theorem and Corollary~\ref{cor.lg-sip}, we can have at least one regular pairing different from $\SIP{x}{y}_+$.

The final statement about ``almost uniqueness'' is immediate from the Rademacher theorem~\cite{HF:69}
stating that on a finite dimensional space $X$ the norm (and each Lipschitz function) is almost
everywhere differentiable. If the norm is differentiable at point $y\in X$, then one has
$\SIP{x}{y}_-=\SIP{x}{y}_+=\WP{x}{y}$ for all $x\in X$ in view of~\eqref{eq.dominated}. This finishes the proof.

\begin{remark}
  Notice that we in fact have proved that a normed space that is not
  Gateaux differentiable admits infinitely many different regular pairings,
  because a convex combination of two regular pairings is, obviously, also
  a regular pairing. \oprocend
\end{remark}

\section{Proof of Lemma~\ref{lem.wp-and-sip}}\label{app.proof-lem-wp-sip}
\textbf{The real space case.}  Implication $\ref{lem1:real}\implies\ref{lem1:hom}$ is straightforward.

To prove the inverse implication, notice first that $\ref{lem1:hom}$ entails the additivity~\eqref{eq.additive} of $\SIP{x}{y}=\WP{x}{y}$ in view of the inequalities
\begin{align*}
  \WP{x_1+x_2}{y}&\leq\WP{x_1}{y}+\WP{x_2}{y}=\WP{x_1}{y}+\WP{x_1+x_2-x_1}{y} \\
  & \leq \WP{x_1}{y}+\WP{-x_1}{y}+\WP{x_1+x_2}{y}=\WP{x_1+x_2}{y}.
\end{align*}
Due to~\eqref{eq.weak-hom1} and~\eqref{eq.weak-hom2},~\eqref{eq.hom1}
and~\eqref{eq.hom2} are also valid.  Finally,~\eqref{eq.posdef-sip}
and~\eqref{eq.cauchy-sip} are ensured by the axioms of WP. This finishes
the proof in the real space case.

\textbf{The complex space case.} Implication $\ref{lem1:real+}\Longrightarrow\ref{lem1:hom+}$ is straightforward.

To prove the inverse implication $\ref{lem1:hom+}\Longrightarrow\ref{lem1:real+}$, notice first that for any LG pairing $\SIP{\cdot}{\cdot}$ one has $\SIP{\jmath x}{y}=\jmath\SIP{x}{y}$.
Hence, if LG pairing satisfying~$\ref{lem1:real+}$ does exist, we have $\WP{x}{y}=\realpart\SIP{x}{y}$ for all $x,y\in X$,
whence $\WP{\jmath x}{y}=-\imagpart\SIP{x}{y}$. Hence, if an LG pairing satisfying~$\ref{lem1:real+}$ exists, it is uniquely determined from
\begin{equation}\label{eq.wp-to-sip}
\SIP{x}{y}:=\WP{x}{y}-\jmath\WP{\jmath x}{y}\quad\forall x,y\in X.
\end{equation}
It remains to check that the conditions $\ref{lem1:hom+}$ indeed make~\eqref{eq.wp-to-sip} an LG pairing.
The proof of additivity~\eqref{eq.additive} retraces the real case, as well as the proof of
~\eqref{eq.hom1} and~\eqref{eq.hom2} for $\lambda\in\real$. Since by definition $\SIP{\jmath x}{y}=\jmath\SIP{x}{y}$ and the pairing~\eqref{eq.wp-to-sip} is additive,~\eqref{eq.hom1} also holds for all $\lambda\in\complex$. To prove~\eqref{eq.hom2}, it suffices to consider $\lambda$ such that $|\lambda|=1$. Using~\emph{\ref{lem1:hom+}}, one has $\WP{x}{\lambda y}=\WP{\bar\lambda x}{\lambda\bar\lambda y}=\WP{\bar\lambda x}{y}$ for all $x,y\in X$. Hence
\begin{equation*}
\SIP{x}{\lambda y}=\WP{\bar\lambda x}{y}-\jmath\WP{\jmath\bar\lambda x}{y}=\SIP{\bar\lambda x}{y}\overset{\eqref{eq.hom1}}{=}\bar\lambda\SIP{x}{y}.
\end{equation*}
Since $\SIP{x}{x}=\WP{x}{x}$, ~\eqref{eq.posdef-sip} and~\eqref{eq.cauchy-sip} also hold. We have proved that~\eqref{eq.wp-to-sip} is the unique LG pairing,
satisfying~$\ref{lem1:real+}$.

\section{Proof of Lemma~\ref{lem.max-pairing-jmt}}\label{app.proof-lem-max-pairing}
For $y=\vect{0}$, the statement is obvious. Hence, we assume that $\|y\|_{\infty}>0$ without loss of generality.

Notice that if $x_i,y_i$ change their signs simultaneously for some $i$, neither norms $\|y+hx\|_{\infty}$, $h\in\mathbb{R}$, nor the product $y_ix_i$ change.
In other words, both the max pairing and the upper JMT pairing for the $\ell_{\infty}$ norm are invariant to transformations
$x\mapsto Ex$, $y\mapsto Ey$, where $E$ is a diagonal matrix with diagonal entries equal to $\pm 1$.

Hence, it suffices to prove the lemma for the case where the second
argument is a nonnegative vector: $y_i\geq 0$ for each $i$. In this
situation, $I_{\infty}(y)=\{i:y_i=\|y\|_{\infty}=\max_i y_i\}$;
furthermore, it can be shown that $\|y+hx\|_{\infty}=\max_i(y_i+hx_i)$ for
$h>0$ being sufficiently small\footnote{Indeed, it is easily seen that
$I_{\infty}(y+hx)\subseteq I_{\infty}(y)$ for $h$ being small. Recalling
that $y\ne\vect{0}$ and thus $y_i>0$ for $i\in I_{\infty}(y)$, one has
$y_i+hx_i>0$ for all $i\in I_{\infty}(y+hx)$ when $h$ is small enough.}.
Hence, the upper JMT pairing is found as
\[
[x,y]_{\infty,+}=\|y\|_{\infty}\lim_{h\to 0+}\frac{f(y+hx)-f(y)}{h}=\|y\|_{\infty}\frac{\partial f(y)}{\partial x},
\]
where $f(y):=\max_{i\in\{1,\ldots,n\}}f_i(y)$ and $f_i(y):=y_i$. Using the Danskin theorem\footnote{Formally, the Danskin theorem in~\cite{JMD:66} is formulated for the minimum, but its extension to the maximum is straightforward.}~\cite[Theorem~I]{JMD:66}, the latter directional derivative can be computed as follows
\[
\frac{\partial f}{\partial x}=\max_{i\in I_{\infty}(y)}\frac{\partial f_i}{\partial x}=\max_{i\in I_{\infty}(y)}x_i.
\]
Here we use the fact that $I_{\infty}(y)=\{i:f_i(y)=y_i=f(y)=\|y\|_{\infty}\}$.
Hence, assuming that $y$ is a nonnegative vector with $\|y\|_{\infty}>0$, for each $x$ we have
\[
[x,y]_{\infty,+}=\|y\|_{\infty}\frac{\partial\|y\|}{\partial x}=\|y\|_{\infty}\max_{i\in I_{\infty}(y)}x_i=\max_{i\in I_{\infty}(y)}x_iy_i=\WP{x}{y}_{\infty},
\]
which, in view of the previous remarks, finishes the proof.

\section{Proof of Lemma~\ref{lem:polyh-max-pairing}}\label{app.proof-lem-polyh-max-pairing}
	It is immediate to verify that the norm associated to $\WP{\cdot}{\cdot}_W$ is $\norm{\cdot}{W}$.
	The upper JMT pairing for this norm is found via Lemma~\ref{lem.max-pairing-jmt} as follows
	\[
	\begin{aligned}
		[x,y]_{W,+}&=\|y\|_W\lim_{h\to 0+}\frac{\|y+hx\|_W-\|y\|_W}{h}=\|Wy\|_{\infty}\lim_{h\to 0+}\frac{\|W(y+hx)\|_{\infty}-\|Wy\|_{\infty}}{h}\\
		&=[Wx,Wy]_{\infty,+}\overset{\eqref{eq.max-pairing-jmt}}{=}\WP{Wx}{Wy}_{\infty}=\WP{x}{y}_W.
	\end{aligned}
	\]
	This concludes our proof of the first statement: the polyhedral max pairing is the upper JMT pairing (and thus is regular).
	
	To establish the equality~\eqref{opt-problem-lognorm}, we first construct
	a matrix $H^*$ that satisfies~\eqref{eq:weird-comm} and the inequality
	$\mu_{\infty}(H^*)\leq\mu_W(A)$. This construction is broadly known in
	the literature~\cite{APM-ESP:86,AB-CB:88,EBC-JCH:93} and exploits the
	homogeneous Farkas lemma~\cite[Lemma~1.4.1]{ABT-AN:01}.
	
	Denote the row vectors of $W$ by $w_1^{\top},\ldots,w_m^{\top}$ and let
	$b_0:=\mu_{\infty}(A)$. Choose an index $i\in\{1,\ldots,m\}$ and consider
	a vector $x\in\mathbb{R}^n$ obeying $2(m-1)$ inequalities
	\begin{subequations}\label{eq.aux-ineq}
		\begin{gather}
			w_i^{\top}x\geq w_j^{\top}x\quad\forall j\ne i, \label{eq.aux-ineq1}\\
			w_i^{\top}x\geq -w_j^{\top}x\quad\forall j\ne i.\label{eq.aux-ineq2}
		\end{gather}
	\end{subequations}
	Then, one obviously has $\|x\|_W=\|Wx\|_{\infty}=w_i^{\top}x$. Using the
	definition of max pairing and the Lumer inequality, we compute
	\[
	b_0(w_i^{\top}x)^2=b_0\|x\|_W^2\overset{\eqref{eq.lumer1}}{\geq}\WP{Ax}{x}_W=\WP{WAx}{Wx}_{\infty}\overset{\eqref{def:WP-1+infty}}{\geq}
	(w_i^{\top}Ax)(w_i^{\top}x),\] which, dividing by\footnote{As we have
		proved, if $\|x\|_W=\|Wx\|_{\infty}=0$, then $x=0$
		and~\eqref{eq.aux-ineq3} retains its validity.}
	$w_i^{\top}x=\|x\|_W$, leads to the following inequality:
	\begin{equation}\label{eq.aux-ineq3}
		(b_0w_i^{\top}-w_i^{\top}A)x\geq 0.
	\end{equation}
	In other words, the system of inequalities~\eqref{eq.aux-ineq} implies
	the inequality~\eqref{eq.aux-ineq3} for each $i\in\until{m}$. The
	homogeneous Farkas lemma~\cite[Lemma~1.4.1]{ABT-AN:01} thus guarantees
	the existence of numbers $\lambda_{ij}^+,\lambda_{ij}^-\geq 0$ (defined
	for $1\leq i,j\leq m$ and $i\ne j$) such that
	\begin{equation*}
		b_0w_i^{\top}-w_i^{\top}A =
		\sum\nolimits_{j\ne i}\lambda_{ij}^+(w_i^{\top}-w_j^{\top})+\sum\nolimits_{j\ne i}\lambda_{ij}^-(w_i^{\top}+w_j^{\top}),
	\end{equation*}
	which is equivalently written as follows:
	\begin{equation}\label{eq.aux3}
		\Bigl(b_0-\sum\nolimits_{j\ne i}\left(\lambda_{ij}^++\lambda_{ij}^-\right)\Bigr) w_i^{\top}
		+\sum\nolimits_{j\ne i}\left(\lambda_{ij}^+-\lambda_{ij}^-\right)w_j^{\top}=w_i^{\top}A.
	\end{equation}
	The desired matrix $H^*$ can now be defined by
	\[
	H^*_{ij}=
	\begin{cases}
		\lambda_{ij}^+-\lambda_{ij}^-,\quad &i\ne j,\\
		b_0-\sum_{j\ne i}(\lambda_{ij}^++\lambda_{ij}^-),\quad &i=j
	\end{cases}
	\]
	Recalling that $\lambda_{ij}^{\pm}\geq 0$ and using
	Table~\ref{table:equivalences}, it is easily shown that
	\[
	H^*_{ii}+\sum_{j\ne i}|H^*_{ij}|\leq b_0\quad\forall i\in\until{m},
	\]
	and hence $\mu_{\infty}(H^*)\leq b_0=\mu_W(A)$. Note that the $i$th row of matrix $H^*W$ is found as
	\[
	(H^*W)_{i\cdot}=H^*_{ii}w_i^{\top}+\sum_{j\ne i}H_{ij}^*w_{j}^{\top}\overset{\eqref{eq.aux3}}{=}w_i^{\top}A,
	\]
	for all $i=1,\ldots, m$, which is nothing else than the $i$th row of $WA$. Hence, $H^*$ satisfies~\eqref{eq:weird-comm} and the
	inequality $\mu_{\infty}(H^*)\leq\mu_W(A)$.
	
	Next, we prove the converse inequality $b_0=\mu_W(A)\leq\mu_\infty(H)$
	for all $H$ satisfying~\eqref{eq:weird-comm}.  To do so, we proceed as
	follows. For each $x\in\real^n\setminus\{0\}$, we compute
	\begin{align*}
		\WP{Ax}{x}_W &= \WP{WAx}{Wx}_\infty \overset{\eqref{eq:weird-comm}}{=} \WP{\barAA Wx}{Wx}_\infty	
		\overset{\text{Lumer's property~\eqref{eq.lumer-a}}}{\leq}
		\mu_\infty(H) \norm{Wx}{\infty}^2 \\
		& = \mu_\infty(H) \norm{x}{W}^2 .
	\end{align*}
	By taking the maximum over $x$ such that $\norm{x}{W}=1$, Lumer's
	property for the regular pairing $\WP{\cdot}{\cdot}_W$ results in
	$\mu_W(A)=\sup\nolimits_{\norm{x}{W}=1}\WP{Ax}{x}_W \leq \mu_\infty(H)$.
	
	We have proved that $b_0=\mu_W(A)\leq\mu_\infty(H)$ whenever
	$H$ satisfies~\eqref{eq:weird-comm}. In particular,
	$b_0\leq\mu(H^*)$. This means that the matrix $H^*$ constructed in the
	first part of the proof is a minimizer in the
	LP~\eqref{opt-problem-lognorm}, and the minimum is
	$\mu_{\infty}(H^*)=b_0=\mu_W(A)$.

\end{document}